\documentclass[12pt]{article}
\usepackage{amsmath}
\usepackage{amsthm}
\usepackage{amssymb}
\usepackage{eucal}
\usepackage{graphicx,graphics}
\numberwithin{equation}{section}
\usepackage{graphicx,graphics}
\usepackage{pdfsync}
\def \dis {\displaystyle}

\def \into {\int_\Omega}
\def \confai {-\kern -.5em\rightharpoonup}
\def \cqfd {\hfill$\Box$}

\def \ep {\varepsilon}
\def \om {\omega}
\def \Om {\Omega}

\def \ZZ {\mathbb Z}
\def \RR {\mathbb R}
\def \CC {\mathbb C}
\def \beq {\begin{equation}}
\def \eeq {\end{equation}}
\def \ba {\begin{array}}
\def \ea {\end{array}}
\def \bs {\bigskip}

\def \ecart {\noalign{\medskip}}
\def\intm {-\hskip-6pt-\hskip-16pt\int}
\newtheorem{Thm}{Theorem}[section]

\newtheorem{Pro}[Thm]{Proposition}
\newtheorem{Lem}[Thm]{Lemma}

\newtheorem{Adef}[Thm]{Definition}

\newtheorem{Arem}[Thm]{Remark}
\newenvironment{Rem}{\begin{Arem}\it}{\end{Arem}}
\newtheorem{Aexa}[Thm]{Example}

\newtheorem{Anot}[Thm]{Notation}

\def \refe #1.{(\ref{#1})}
\def \reff #1.{figure~\ref{#1}}
\def \refs #1.{Section~\ref{#1}}
\def \refss #1.{Subsection~\ref{#1}}
\def \refD #1.{Definition~\ref{#1}}
\def \refT #1.{Theorem~\ref{#1}}
\def \refL #1.{Lemma~\ref{#1}}
\def \refC #1.{Corollary~\ref{#1}}
\def \refP #1.{Proposition~\ref{#1}}
\def \refPt #1.{Properties~\ref{#1}}
\def \refR #1.{Remark~\ref{#1}}
\def \refE #1.{Example~\ref{#1}}
\def \refN #1.{Notation~\ref{#1}}
\title{A viscoelastic model for skin
\\
via homogenization theory.}
\author{
\footnotesize
\centerline{\begin{tabular}{cc}
\Large Juan Casado-D\'iaz
\\
 Dpto. de Ecuaciones Diferenciales y An\'alisis Num\'erico
\\
 Universidad de Sevilla
\\
 jcasadod@us.es
\end{tabular}}
}

\begin{document}

\maketitle
\begin{abstract}
We carry out the homogenization of a fluid-structure interaction problem consisting in 
the periodic inclusions of a viscous fluid in an elastic body. We get a macrostructure model where the body behaves as a viscoelastic material with a long-range memory term. Our aim is not only to get this limit problem but also to study its main properties. Using the micro-structure variables
it is simple to check that it satisfies an energy conservation law assuring in particular the existence and uniqueness of solution. The difficulty is to characterize these properties using only the macroscopic system. We prove that the nonlocal term is given through  a convolution  kernel which exponentially decreases to zero and satisfies some positive conditions which we write in terms of the Laplace transform. These conditions can be used to directly prove the existence and uniqueness of solution. The results apply to modeling the mechanical behavior of skin as an indicator of what kind of models we should use. In a naive interpretation the fluid inclusions represent the cells and the elastic medium the extracellular matrix. 
\end{abstract}

\par\bs\noindent

\section{Introduction.}
In the modeling of the mechanical properties of skin, it is usual to considered that it behaves as a viscoelastic material, i.e. an elastic material where the stress tensor contains a memory term. The most simple case corresponds to  an instantaneous memory term in such a way that the elastic deformation $u$ solves a local partial differential system of the form
$$\rho\partial^2_{tt}u-{\rm div}_x\Big(Ae_x(u)+Be_x(\partial_tu)\Big)=f,$$
with $\rho>0$ the density of the material, $f$ an exterior force, $e_x(u)$ the symmetric part of the spatial derivative of $u$ (strain tensor) and $A,B$ symmetric positive tensors on the space of symmetric matrices.  The instantaneous memory term corresponds to $-{\rm div}_x(Be_x(\partial_tu))$. It was noted in \cite{Mur}, chapter 10, that it may be better to work with a long-range memory term. In the most simple case this corresponds to
$$\rho\partial^2_{tt}u-{\rm div}_x\left(Ae_x(u)+\int_0^tS(t-s)e_x(u)ds\right)=f.$$
Observe that integrating by parts the integral term in this equation and assuming that $S$ is quickly decreasing in $s$, we get that the behavior of the solutions of this second equation approximates that of the previous one.\par 
Our purpose in the present paper is to prove how this kind of memory terms can be obtained by using inclusions of a viscous fluid in a material governed by  the usual linear elastic system. Namely, we consider a very simple and naive model for skin deformation, where the cells are represented by small connected domains filled of a viscous fluid and  periodically distributed. They are surrounded by an elastic material, which represents what in Biology is known as the extra-celular matrix. The behavior of the fluid inside the cells is modeled by the Stokes system, while the deformation of the extra-celular matrix is given by the linear elasticity system. More exactly, recalling that the elasticity system is written in Lagrange coordinate, the fluid-structure interaction model that we consider consists in writing  also the Navier-Stokes system in Lagrange coordinates and then linearizing the resulting equations. This provides a partial differential system in a bounded reference domain $\Om\subset \RR^N$, $N\geq 2$ ($N=2,3$ in the applications), where the cells ocupe a family of connected open sets periodically distributed and  the unknowns are the elastic deformation of the extra-celular matrix and the velocity  and pressure of the fluid inside the cells. The system is accomplished with two boundary conditions on the interface corresponding to the continuity of the velocity and the equilibrium of internal forces. See (\ref{sistho}) below. The system is hyperbolic in the solid part and parabolic in the fluid one. 
We assume that the period is given by a small parameter $\varepsilon>0$ and that the set of cells is given as 
$$\Om_\ep^l:=\bigcup_{\ep k+\ep Y\subset\Om\atop k\in\ZZ^N} (\ep k+\ep O),$$
with $Y$ the unitary cube in $\RR^N$ and $O$ a connected open set  sufficiently smooth and such that $\overline O\subset Y$. We also denote by $\Om_\ep^s$ the solid part. Our purpose is to pass to the limit when $\ep$ tends to zero in this partial differential system to get a macroscopic model.\par\smallskip
Taking initial conditions  bounded in the good space to have a unique solution, and extracting a convergent subsequence  (see Proposition \ref{Prolic}),
we prove in Theorem \ref{Threpvel}
that denoting by $u_\ep$ the elastic deformation  of the of the solid part and by $v_\ep$ the velocity of the fluid, we have
 $$u_\ep\stackrel{\ast}\rightharpoonup u\ \hbox{ in }W^{1,\infty}(0,T;L^2(\Om))^N\cap L^\infty(0,T;H^1_0(\Om))^N,\quad {1\over |O|}v_\ep\chi_{\Om_\ep^l}\stackrel{\ast}\rightharpoonup \partial_t u\ \hbox{ in }L^\infty(0,T;L^2(\Om))^N,$$
 where $u$ solves a partial differential system of the form
\beq\label{pbliI}\ba{l}\dis \rho^{eff}\partial^2_{tt}u-{\rm div}_x\Big(A_h e_x(u)+\int_0^t S(t-s)e_x(u)\,ds\\ \ecart\dis
\qquad\qquad\qquad\qquad\qquad +\int_{\partial O} R(t,y)\big(e_x(u^0)y+u^1_0\big)\,d\sigma(y)\Big)=f\ \hbox{ in }(0,T)\times \Om.\ea\eeq
In this system $f$ is the exterior force, $\rho^{eff}$ is the mean value of the solid and fluid densities, $u^0$ the initial condition for $u$, and $u_1^0$ a function in $H^{1\over 2}(\partial O)$ depending also on the initial conditions. The integral on $\partial O$ must be understood as the duality product of $R(t,.)\in H^{-{1\over 2}}(\partial O)^{N\times N\times N}$ and $e_x(u^0)y+u^1_0\in H^{1\over 2}(\partial O)^N$.  The proof uses the unfolding method (\cite{ArDoHo}, \cite{Cas}, \cite{CiDaGr}), which is strongly related with the two-scale convergence method in periodic homogenization (\cite{All}, \cite{Ngu}).\par
As we said above the microscopic model considered  is very simple from the biological point of view, but we are not only interested in getting this homogenized system (\ref{pbliI}). Our aim is to impose the properties of this limit problem to general models describing the mechanical properties of skin. \par
We prove  that the functions $R$ and $S$ in (\ref{pbliI}) have an exponential decay at infinity. This means that although the model has memory, only recent evens have influence. We also prove that the tensor $A_h$ is symmetric and positive and $S$ is symmetric. Finally, we deduce that 
the limit model  satisfies an energy conservation law (see (\ref{energS})) depending on the microscopic structure, where
due to the viscosity of the fluid, the energy decreases with time. The difficulty is to write this condition using only the macroscopic model. We prove that this provides
 a positivity condition for the sum of $A_h$ and the memory term which we write in terms of the Laplace transform of $S$ (see (\ref{proLS1}), (\ref{proLS2})).
 As a consequence of  this conditions we deduce the existence of $c>0$, such that $A_h$ and $S$ satisfy the inequality
$$\ba{l}\dis\int_0^T \Big(A_hE+\int_0^tS(t-s)E(s)ds\Big){dE\over dt}\,dt\\ \ecart\dis\geq c|E(T)|^2+c\int_0^T \Big|\int_0^te^{-(t-s)}{dE\over ds}ds\Big|^2dt+c\Big|\int_0^Te^{-(T-s)}{dE\over ds}ds\Big|^2,\ea$$
for every $E\in W^{1,\infty}(0,\infty;\CC^N_s)$, such that $E(0)=0$. It allows us to prove the existence and uniqueness of solution for (\ref{pbliI}).
 \par
Although in our knowledge this is the first work where a memory term is rigorously obtained from the microscopic mixture of fluids and solids, the apparition of  these terms in homogenization has before been obtained in other situations. In this way, it was proved in \cite{San}, chapter 6, that the limit of the solutions of a viscoelastic problem with an instantaneous memory term where the coefficients are rapidly oscillating gives in general a viscoelastic  problem with a long-range memory term with exponential decay at infinity.
In \cite{ArDoHo}, it is considered a geometry similar to that of this paper, but instead of taking an interaction fluid-structure problem,  the authors consider the heat equation where the diffusion coefficients are of order $\ep^2$ in the periodic inclusions  and of order one in the rest. The corresponding limit equation contains a zero-order memory term which exponentially decreases to zero. If instead of considering the heat equation we consider the linear elasticity system, then we get a viscoelastic behavior in the limit with a stress tensor containing a memory term (\cite{BrCa}, see also \cite{AGMR} for the case of the elasticity system with a fixed frequency). Since the elasticity system conserves the energy, this memory term does not converges exponentially to zero as in our case. Another kind of equations which provide a memory term corresponds to the wave equation (or the elasticity system) with a highly oscillating drift term (\cite{BrCa2}, \cite{CCMM}, \cite{DMSa}). In this case we get a nonlocal drift term depending on a dependence cone.\par\medskip
The paper is organized as follows:\par\medskip
In Section \ref{not} we introduce some notation that we will use along the paper.\par\medskip
In Section \ref{seFlSt}, we present the fluid-structure interaction model and we recall an  existence and uniqueness theorem for such model. It will provide the a priori estimates that we will need to carry out the homogenization process.\par\medskip
In Section \ref{secvimo}, we present the periodic model mentioned above and we state and prove the main results of the paper described above.\par
\section{Notation.} \label{not}
\begin{itemize}
\item We denote by $\RR^{N\times N}_s$, $\CC^{N\times N}_s$ the spaces of symmetric matrices of dimension $N\times N$ with coefficients in $\RR$ and $\CC$ respectively.
\item The transposed of a matrix $\xi$  is denoted by $\xi^T$.
\item For two vectors $a,b\in \RR^N$ or $\CC^N$ we use the notation $a\cdot b$ to mean
$\sum_{i=1}^N a_ib_i$,
with $a_i,b_i$ the components of $a$ and $b$. For vectors in $\RR^N$ this is just the usual scalar product of $a$ and $b$. In the complex case, such scalar product is given by $a\cdot\overline{b}$ with $\overline b$ the conjugate of $b$.  Analogously, for two matrices $\xi,\eta\in \RR^{N\times N}$  or $\CC^{N\times N}_s$ we write $\xi:\eta:={\rm trace}(\xi\eta^T)=\sum_{i,j=1}^N\xi_{ij}\eta_{ij}$.
\item For a function $u$, we denote by $e(u)$ the symmetric part of the derivative of $u$, i.e. $e(u)=(Du+Du^T)/2$.
\item For a measurable set $U\subset \RR^N$, we denote by $\chi_U$ the characteristic function of $U$. For a measurable function $u$ defined in $U$, we denote by $u\chi_U$ the extension by zero of $u$ to $\RR^N$.\item For two normed spaces $X_1,X_2$, we denote by ${\cal L}(X_1;X_2)$ the space of continuous linear functions from $X_1$ into $X_2$. When $X_1=X_2$, we just write ${\cal L}(X_1)$.
\item We denote by $e_1,...,e_N$ the elements of the canonical basis in $\RR^N$.
\item For a normed space $X$, we denote by $I$ the identity operator from $X$ into $X$. The space $X$ will be understood by the context.
\item For an open set $U\subset\RR^N$ whose boundary is locally the graph of a Lipschitz function, and a subset $\Gamma\subset\partial U$, we define $H^1_\Gamma(U)$ as the space of functions in $H^1(U)$ which vanish on $\Gamma$. When $\Gamma=\partial U$, we write $H^1_0(U)$. 
\item For $g\in H^{-{1\over 2}}(\Gamma)$ and $u\in H^{1\over 2}(\Gamma)$, we use the notation
$$\int_{\Gamma} gu\,d\sigma$$
to mean de duality product $\langle g,u\rangle_{H^{-{1\over 2}}(\Gamma),H^{{1\over 2}}(\Gamma)}$.
Remark that in general $g$ is not a function and therefore this duality product is not well defined through the previous integral.
\item For a Banach space $X$ over $\CC$ and $f\in L^1(0.T;X)$ for every $T>0$ such that there exists $\gamma\in\RR$ with $e^{- \gamma s}f\in L^1(0,\infty;X)$, we recall the existence of $\alpha\geq\infty$ such that $e^{-a s}f\in L^1(0,\infty;X)$ if $a>\alpha$ and $e^{-a s}f\not\in L^1(0,\infty;X)$ if $a<\alpha$.
In these conditions, we define the Laplace transform of $f$, $Lf$ by
$$Lf(z)=\int_0^\infty e^{-zt}f(t)\,dt,\quad \forall\, z\in \CC,\ \hbox{ with }Re(z)>\alpha.$$
\item For $U\subset\RR^N$ open, we denote by ${\cal D}(U)$ the space of $C^\infty$-functions with compact support in $U$. We recall that it is endowed with the inductive limit topology of the spaces ${\cal D}_K(U)$ of functions in ${\cal D}(U)$ with support in $K\subset U$ compact. The topological dual of ${\cal D}(U)$ is the space of distributions in $U$, and it is represented by ${\cal D}'(U)$.
\item We denote by $C$ a generic positive constant which can change from line to line.
\item We denote by $O_\ep$ a generic sequence which tends to zero with $\ep$ and can change from line to line.
\item We denote by $Y$ the unit cube $(-1/2,1/2)^N$ in $\RR^N$. We use the index $\sharp$ to mean periodicity with respect to $Y$. For example, $H^1_\sharp(Y)$ is the space of functions $u$ in the Sobolev space $H^1_{loc}(\RR^N)$ such that $u(y+k)=u(y)$ for a.e. $y\in Y$ and every $k\in\ZZ$. More generally, we use this symbol for functions defined in a subset of $Y$. For example, for an open subset $U\subset Y$, $H^1_\sharp(U)$ is the space of functions $u$ in 
$H^1_{loc}(\cup_{k\in\ZZ}(k+U))$ such that $u(y+k)=u(y)$ for a.e. $y\in U$ and every $k\in\ZZ^N$.  When $u$ satisfies $u(x+e_i)=-u(x)$, for every $i\in \{1,...,N\}$ and a.e. $x\in U$ we say that $u$ is anti-periodic in $Y$.
\item For $k\in\ZZ^N$ and $\ep>0$, we denote $C_k^\ep:=\ep (k+Y)$.
\item We define $\kappa:\RR^N\to \ZZ^N$ by
$$\kappa(x)\in \ZZ^N,\quad x-\kappa(x)\in Y,\ \hbox{ a.e. }x\in\RR^N,$$
i.e. $\kappa(x)$ gives the center of the cube $C_k^1$ containing $x$. Observe that more generally, for $\ep>0$, $\ep\kappa(x/\ep)$ is the center of the cube $C_k^\ep$ containing $x$.
\end{itemize}
\section{ A fluid-structure interaction model.}
\label{seFlSt}
We present in this section the partial differential system corresponding to the simple fluid-structure interaction model that we will use in the following one. We  also recall some results about the existence and uniqueness of solution for this problem. \par
We consider two open sets $\om,\Om\subset \RR^N$  such that $\overline\om\subset\Om$ and such that their respective boundaries can be locally written as the graph of a Lipschitz function. We assume that $\om$ is filled with a viscous fluid governed by the Navier-Stokes system, while $\Om\setminus\overline \om$ is the solid part, governed by the linear elasticity system.\par\par
 We recall that the elasticity system is usually written in Lagrange coordinates, i.e. the coordinates refer to a fixed  geometry in which there is no internal forces. In our case this geometry is given by $\Om\setminus\overline\om$. The linear elasticity system is then given by
 \beq\label{elsi}\rho^s\partial^2_{tt}u-{\rm div}_x (Ae_x(u))=f\ \hbox{ in }(0,T)\times(\Om\setminus\overline\om),\eeq
 where $T>0$ is the final time, $\rho^s>0$  the density of the solid part,  $A\in{\cal L}(\RR^{N\times N}_s)$, the stress tensor, is an element of ${\cal L}(\RR^{N\times N}_s)$ which we assume to be positive and symmetric, i.e.  such that
\beq\label{proeA}\exists\,\alpha>0,\ \  A\xi:\xi\geq \alpha|\xi|^2,\ \forall\,\xi\in\RR^{N\times N}_s,\qquad A\xi:\eta=A\eta:\xi,\ \ \forall\,\xi,\eta\in \RR^{N\times N}_s,\eeq
and $f$, the exterior force, is a function in $(0,T)\times\Om$.\par
We recall that a particular, and important,  case corresponds to the isotropic elastic materials where  there exist $\lambda^s,\mu^s>0$, the Lam\'e's constants, such that the stress tensor $A$ is given by
\beq\label{isoma} A\xi=\lambda^s{\rm trace}(\xi)I+2\mu^s \xi,\quad \forall\, \xi\in \RR^{N\times N}_s.\eeq
For our purposes we will not need $A$ to be isotropic.\par
The function $u=u(t,x)$ in (\ref{elsi}) represents the deformation. This means that the particle that was in position $x$ at the initial time moved to  $x+u(t,x)$ at time $t$. \par
The movement of the fluid is usually written in Euler coordinates. Thus, these equations are not given in the fixed open set $\om$ but in an open set $\om(t)$ which variess with $t$, and are given by the Navier-Stokes system
\beq\label{EcfLc}\left\{\ba{l}\dis \rho^l\big(\partial_tv+(v\cdot\nabla_y)v\big)-2\mu\,{\rm div}_ye_y(u)+\nabla_y p=F\ \hbox{ in }\big\{(t,x): t\in (0,T),\ y\in \om(t)\big\}\\ \ecart\dis
{\rm div}_y\, v=0\ \hbox{ in }\big\{(t,x): t\in (0,T),\ y\in \om(t)\big\},\ea\right.\eeq
 where $\rho^l>0$ is the density of the fluid, $\mu>0$ the viscosity of the fluid and $F=F(t,y)$ the exterior force. In these variables $v(t,y)$ is the velocity of the particle that at time $t$ occupies the position $y$. This position varies with time contrarily to the Lagrange coordinates. The function $p$ represents the pressure of the fluid.\par
  The Lagrange coordinates $(t,x)$ and the Euler coordinates $(t,y)$ are  linked by 
 $$y=x+u^l(t,x),$$ 
 with $u^l$ the deformation of the fluid in Lagrange coordinates, i.e. the solution of
 $$\partial_tu^l(t,x)=v(t, x+u^l(t,x)),\quad u^l(0,x)=x.$$
 In the following, we consider  a simple fluid-structure interaction model. It consists in writing the Navier-Stokes system  in Lagrange coordinates and then to carry out a linearization around the null function. This means that we are assuming  the deformations to be small. The boundary conditions at the interface correspond to the the equilibrium of forces and the continuity of the deformation between the fluid and the solid parts. Adding the corresponding initial conditions and assuming to simplify that the solid part is fixed at the boundary of $\Om$, we get the system  
 \beq\label{promo} \left\{\ba{l}\dis \rho^s\partial^2_{tt}u-{\rm div}_x (Ae_x(u))=f\ \hbox{ in }(0,T)\times(\Om\setminus\overline\om)\\ \ecart\dis \rho^l\partial_tv-2\mu\, {\rm div}_x e_x(v)+\nabla_x p=f\ \hbox{ in }(0,T)\times\om\\ \ecart\dis
 {\rm div}_xv=0\ \hbox{ in }(0,T)\times\om\\ \ecart\dis   \partial_t u=v,\ \ 
 Ae_x(u)\nu=(2\mu e_x(v)-pI)\nu\quad \hbox{on }(0,T)\times\partial\om\\ \ecart\dis
 u=0\ \hbox{ on }(0,T)\times\partial\Om
\\ \ecart\dis 
 u_{|t=0}=u^0,\ (\partial_t u)_{|t=0}=u^1\quad\hbox{in }\Om\setminus\overline\om,\qquad v_{|t=0}=v^0\ \hbox{ in }\om,\ea\right.\eeq
 where $\nu$ denotes the outside normal unitary vector on $\partial\om$.
 \begin{Rem}\label{dened} In fluid mechanics it is usual to denote the symmetric part of the derivative of the velocity (the strain rate tensor) by $d_x(v)$ instead of $e_x(v)$, but because from the Mathematical point of view this is the same differential operator that the one used to define the strain tensor $e_x(u)$, we prefer to use the same notation for both of them.
\end{Rem}
In order to get an existence and uniqueness result for (\ref{promo}), let us take $f$, $u^0$, $u^1$ and  $v^0$  satisfying
 \beq\label{supdi}\ba{c}\dis f\in L^1(0,T; L^2(\Om))^N,\ u^0\in H^1_{\partial\Om}(\Om\setminus\overline\om)^N,\ 
 u^1\in L^2(\Om\setminus \overline\om)^N\\ \ecart\dis v^0\in L^2(\om)^N,\ {\rm div}_xv^0=0\hbox{ in }\om.\ea\eeq
We look for solutions  $(u,v,p)$, which satisfy the problem in the distributional sense. Namely, using test functions with null divergence in the fluid part $\om$, we say that $(u,v,p)$ is a solution of (\ref{supdi}) if and only if $(u,v)$ satisfies
\beq\label{pbvsuv} \left\{\ba{l}\dis
u\in W^{1,\infty}(0,T;L^2(\Om\setminus\overline\om))^N\cap L^\infty(0,T;H^1_{\partial \Om}(\Om\setminus\overline\om))^N\\ \ecart\dis  v\in  L^\infty(0,T;L^2(\om))^N\cap L^2(0,T;H^1(\om))^N\\ \ecart\dis
 {\rm div}_xv=0\ \hbox{ in }(0,T)\times\om,\quad
\partial_t u=v\ \hbox{ on }(0,T)\times\partial\om\\ \ecart\dis
u_{|t=0}=u^0,\ (\partial_t u)_{|t=0}=u^1\quad\hbox{in }\Om\setminus\overline\om,\qquad v_{|t=0}=v^0\ \hbox{ in }\om\\ \ecart\dis
{d\over dt}\left(\rho^s\int_{\Om\setminus\overline\om}\hskip-4pt \partial_t u\cdot \varphi\,dx+\rho^l\int_{\om} v\cdot\varphi\,dx\right)\\ \ecart\dis \qquad +\int_{\Om\setminus\overline\om} \hskip-4pt Ae_x(u):e_x(\varphi)\,dx+2\mu\int_\omega e_x(v):e_x(\varphi)\,dx=\into f\cdot\varphi\,dx\ \hbox{ in }{\cal D}'(0,T),\\ \ecart\dis\forall\, \varphi\in H^1_0(\Om)^N,\quad {\rm div}_x\varphi=0\ \hbox{ in }\omega.\ea\right.\eeq
Then, De Rham's theorem proves the existence of a distribution $p\in{\cal D}'((0,T)\times\om)$ such that  the second equation and the second boundary condition in  the fourth line of (\ref{promo}) are satisfied.  \par\medskip
The existence and uniqueness theorem for (\ref{promo}) is given by
 \begin{Thm} \label{teoeu} For every right-hand side $f$ and every initial data $(u^0,u^1,v^0)$ satisfying (\ref{supdi}),  there exists a unique solution $(u,v,p)$ of (\ref{promo}), such that
 \beq\label{regsol} \ba{c}\dis u\in W^{1,\infty}(0,T;L^2(\Om\setminus\overline\om))^N\cap L^\infty(0,T;H^1_{\partial\Om}(\Om\setminus\overline\om))^N\\ \ecart\dis
v\in L^\infty(0,T;L^2(\om))^N\cap L^2(0,T;H^1(\om))^N,\quad p\in W^{-1,\infty}(0,T;L^2(\om)).\ea\eeq
It also satisfies the energy identity
\beq\label{iden}\ba{l}\dis {1\over 2}{d\over dt}\left(\rho^s\int_{\Om\setminus\overline\om}\big(|\partial_tu|^2+Ae_x(u):e_x(u)\big)dx+\rho^l\int_{\om} |v|^2dx\right)\\ \ecart\dis\qquad\qquad+2\mu\int_\om|e_x(v)|^2dx=\int_{\Om\setminus\overline\om}\hskip-4pt f\cdot\partial_t u\,dx+\int_{\om} f\cdot\ v\,dx\ \hbox{ in }(0,T).\ea\eeq
 \end{Thm}\par\medskip\noindent
  {\bf Proof.} Let us use the Hille-Yoshida theorem. Taking $w=\partial_t u$, we write (\ref{promo}) as the first order system
 \beq\label{promo2} \left\{\ba{l}\partial_t u=w\ \hbox{ in }(0,T)\times(\Om\setminus\overline\om)\\ \ecart\dis
  \dis \partial_tw={1\over\rho^s}\,{\rm div}_x Ae_x(u)+{1\over\rho^s}f\ \hbox{ in }(0,T)\times(\Om\setminus\overline\om)\\ \ecart\dis \partial_tv={2\mu\over \rho^l}\, {\rm div}_x e_x(v)-{1\over \rho^l}\nabla_x p+{1\over \rho^l}f\ \hbox{ in }(0,T)\times\om\\ \ecart\dis
 {\rm div}_xv=0\ \hbox{ in }(0,T)\times\om\\ \ecart\dis   w=v,\ \ 
 Ae_x(u)\nu=(2\mu e_x(v)-pI)\nu\quad \hbox{on }(0,T)\times\partial\om\\ \ecart\dis
 u=0\ \hbox{ on }(0,T)\times\partial\Om
\\ \ecart\dis 
 u_{|t=0}=u^0,\ w_{|t=0}=u^1\quad\hbox{in }\Om\setminus\overline\om,\qquad v_{|t=0}=v^0\ \hbox{ in }\om.\ea\right.\eeq
We take ${\cal H}$ as the Hilbert space
 $${\cal H}=\big\{(u,w,v)\in H^1_{\partial\Om}(\Om\setminus\overline\om)^N\times L^2(\Om\setminus\overline\om)^N\times L^2(\om)^N:\quad  {\rm div}_xv=0\ \hbox{ in }\om \},$$
 endowed with the scalar product
 $$\big((u_1,w_1,v_1),(u_2,w_2,v_2)\big)_{\cal H}=\int_{\Om\setminus\overline\om}\Big(Ae_x(u_1):e_x(u_2)+\rho^s w_1\cdot w_2\Big)\,dx+\rho^l\int_\om v_1\cdot v_2\,dx.$$
 Then we define the unbounded operator ${\cal A}:D({\cal A})\subset{\cal H}\to {\cal H}$
 by
 $$\ba{l}\dis D({\cal A}):=\Big\{(u,w,v)\in {\cal H}:\ {\rm div}_x(Ae_x(u))\in L^2(\Om\setminus\overline\om)^N,\ w\in H^1(\Om\setminus\overline\om)^N,\\ \ecart\dis  \exists p\in L^2(\om)\hbox{ with }2\mu\,	{\rm div}_xe_x(u)-\nabla_x p\in L^2(\om)^N,\ w=v,\ Ae_x(u)\nu=(2\mu\, e_x(v)-pI)\nu\ \hbox{ on }\partial\om\Big\},\ea$$
 and (observe that $p$ in the definition of $D({\cal A})$ is unique)
$${\cal A}(u,w,v)=\Big(w,\,{1\over\rho^s}{\rm div}\,Ae(u),\, {2\mu\over \rho^l}\,{\rm div}\,e(v)-{1\over\rho^l}\nabla p\Big),\quad\forall\,(u,w,v)\in D({\cal A}).$$
Taking into account that 
$$\ba{l}\dis \big({\cal A}(u,w,v),(u,w,v)\big)_{\cal H}=\int_{\Om\setminus\overline\om}\big(Ae(w):e(u) +{\rm div}\, (Ae(u))\cdot w\big)dx\\ \ecart\dis +\int_\om \big(2\mu\,{\rm div}\,e(v)-\nabla p\big)\cdot v\,dx
=-2\mu\int_\om |e(v)|^2dx\leq 0,\ea$$
 and that by Lax-Milgram's theorem $I-{\cal A}$ is bijective from $D({\cal A})$ into ${\cal H}$, we can apply Hille-Yoshida theorem to deduce that ${\cal A}$ generates a strongly continuous semigroup ${\cal T}$ in ${\cal H}$ and then that (\ref{promo2}) has a unique solution.\par \medskip
 If $(u^0,u^1,v^0)$ belongs to $D(A)$ and $f$ belongs to $W^{1,1}(0,T;L^2(\Om))^N$ then, the solution $(u,v,p)$ of (\ref{promo}) is such that
 $$u^0\in W^{2,\infty}(0,T;L^2(\Om\setminus\overline\om))^N\cap W^{1,\infty}(0,T;H^1_{\partial\Om}(\Om\setminus\overline\om))^N$$ 
 $$v\in W^{1,\infty}(0,T;L^2(\om))^N\times H^1(0,T;H^1(\om)))^N,\quad p\in L^\infty(0,T;L^2(\om)). $$
Multiplying the first equation in (\ref{promo}) by $\partial_tu$ and integrating in $\Om\setminus\overline\om$, multiplying the second equation by $v$ and integrating in $\om$, and adding both equalities we get (\ref{iden}). The general case  follows  by density. \cqfd
 \begin{Rem} More generally, Theorem \ref{teoeu} holds with the same proof, assuming $\rho^s$, $\rho^l$  in $L^\infty(\Om\setminus\overline\om)$ and $L^\infty(\om)$ respectively, uniformly positive, and $A$ a tensor function in $L^\infty(\Om;{\cal L}(\RR^{N\times N}_s))$, uniformly elliptic.\par 
  The fact that the fluid part $\om$ is surrounded by the solid part
 is not important either, but this is the case that will interested in the following section. \end{Rem} 
 \begin{Rem} Since $v$ belongs to $L^2(0,T;H^1(\om))^N$, it is defined  in $(0,T)\times \partial\om$ as an element of $L^2(0,T;H^{1\over 2}(\partial \om))^N$. Analogously  $u$ in $L^\infty(0,T; H^{1\over 2}(\partial\Om))^N$ implies that $\partial_tu$ belongs to $W^{-1,\infty}(0,T; H^{1\over 2}(\partial\Om))^N$. This gives a sense to the equality $\partial_t u=v$ on $(0,T)\times \partial\om$.\par
 By (\ref{promo}),  we also have that $\partial^2_{tt}u$ belongs to $L^1(0,T;H^{-1}(\Om\setminus\overline\om))^N$ and $\partial_t v$ to $L^1(0,T; W')$ with $W'$ the dual space of
$$W=\big\{w\in H^1_0(\om)^N:\ {\rm div}_xw=0\ \hbox{ in }\om\big\}.$$
This allows us to give a sense to the initial conditions in (\ref{pbvsuv}). Even more, the proof of Theorem \ref{teoeu}  implies that $u$ is in $C^1([0,T];L^2(\Om\setminus\overline\om))^N\cap C^0([0,T];H^1(\Om\setminus\overline\om))^N,$ and $v$ in $C^0([0,T];L^2(\om))^N$. \end{Rem}
As a consequence of (\ref{iden}) we have the following a priori estimates for the solution of (\ref{promo}).
\begin{Pro} \label{estpta} There exists a constant $C>0$ depending  only on $\rho^s$, $\rho^l$, $\alpha$, $|A|$, $\om$ and $\Om$, such that
\beq\label{estpri1}\ba{l}\dis \|u\|_{W^{1,\infty}(0,T;L^2(\Om\setminus\overline\om))^N\cap L^\infty(0,T;H^1_{\partial\Om}(\Om\setminus\overline\om))^N}+\|v\|_{L^\infty(0,T;L^2(\om))^N\cap L^2(0,T;H^1(\om))^N}\\ \ecart\dis +\Big\|\int_0^tp(s,x)ds\Big\|_{L^\infty(0,T;L^2(\om))}\\ \ecart\dis \leq C\Big(\|e_x(u^0)\|_{L^2(\Om\setminus\overline\om)^{N\times N}}+\|u^1\|_{L^2(\Om\setminus\overline\om)^N}+\|v^0\|_{L^2(\om)^N}+\|f\|_{L^1(0,T;L^2(\Om))^N}\Big).\ea\eeq
\end{Pro}
\par\noindent
{\bf Proof.} Integrating in time in  (\ref{iden}) and taking into account Korn's inequality, we deduce the existence of $C>0$ depending  only on $\rho^s$, $\rho^l$, $\alpha$, $|A|$, $\om$ and $\Om$, such that
\beq\label{estpri1}\ba{l}\dis \|u\|_{W^{1,\infty}(0,T;L^2(\Om\setminus\overline\om))^N\cap L^\infty(0,T;H^1_{\partial\Om}(\Om\setminus\overline\om))^N}+\|v\|_{L^\infty(0,T;L^2(\om))^N\cap L^2(0,T;H^1(\om))^N}\\ \ecart\dis \leq C\Big(\|e_x(u^0)\|_{L^2(\Om\setminus\overline\om)^{N\times N}}+\|u^1\|_{L^2(\Om\setminus\overline\om)^N}+\|v^0\|_{L^2(\om)^N}+\|f\|_{L^1(0,T;L^2(\Om))^N}\Big).\ea\eeq
Integrating in time the second equality in (\ref{promo}) we also have
\beq\label{cosp}\rho^l\big(v-v^0)-2\mu\, {\rm div}_xe_x\Big(\int_0^tv(s,x)\,ds\Big)+\nabla_x\int_0^tp(s,x)\,ds=\int_0^tf(s,x)ds\ \hbox{ in }(0,T)\times\om.\eeq
This proves that
$$\nabla_x \int_0^tp(s,x)\,ds\in L^\infty(0,T;H^{-1}(\Om))^N,$$
and then that
\beq\label{regp} p\in W^{-1,\infty}(0,T;L^2(\Om)).\eeq
Moreover, taking into account the smoothness of $\om$, there exists $w\in L^\infty(0,T;H^1_0(\Om))^N$ such that
\beq\label{propwi}{\rm div}_x w=\int_0^tp(s,x)ds\hbox{ in }(0,T)\times\om,\ \ \|w\|_{L^\infty(0,T;H^1_0(\Om))^N}\leq C\Big\|\int_0^tp(s,x)ds\Big\|_{L^\infty(0,T;L^2(\om))},\eeq
where  $C>0$ only depends on $\om$ and $\Om$.
Then, multiplying (\ref{cosp}) by $w$, integrating by parts and using (\ref{promo}), we get
$$\ba{l}\dis\int_\omega\Big|\int_0^tp(s,x)ds\Big|^2dx=\rho^l\int_{\om} \big(v-v^0\big)\cdot w\,dx+2\mu\int_\omega e\Big(\int_0^t v(s,x)ds\Big):e(w)\,dx\\ \ecart\dis+\rho^s\int_{\Om\setminus\overline\om}(\partial_t u-u^1)\cdot w\,dx+\int_{\Om\setminus\overline\om}Ae_x\Big(\int_0^t u(s,x)ds\Big):e_x(w)\,dx-\into\int_0^tf(s,x)\,ds\cdot w\,dx.\ea$$
By (\ref{estpri1}) and (\ref{propwi}), this proves the existence of $C>0$ depending only on $\rho^s$, $\rho^l$, $\alpha$, $|A|$, $\om$ and $\Om$ such that
\beq\label{estpri2}\ba{l}\dis \Big\|\int_0^tp(s,x)ds\Big\|_{L^\infty(0,T;L^2(\om))}\\ \ecart\dis \leq  C\Big(\|e(u^0)\|_{L^2(\Om\setminus\overline\om)^{N\times N}}+\|u^1\|_{L^2(\Om\setminus\overline\om)^N}+\|v^0\|_{L^2(\om)^N}+\|f\|_{L^1(0,T; L^2(\Om))^N}\Big).\ea\eeq
This finishes the proof. \cqfd 
\begin{Rem} \label{reep} If $p$ is a distribution in $W^{-1,\infty}(0,T;L^2(\om))$ then, it has no sense the integral
$$\int_0^tp(s,x)\,ds.$$
To be rigorous, (\ref{cosp}) must be understood in the following way  
$$\int_0^tf(s,x)ds-\rho^l\big(v-v^0)+2\mu\, {\rm div}_xe_x\Big(\int_0^tv(s,x)\,ds\Big)$$
is a distribution in $L^\infty(0,T;H^{-1}(\om))^N$ which vanishes on the functions in $L^1(0,T;H^1_0(\om))^N$ with null divergence. Therefore, by De Rham's theorem, there exists $P\in L^\infty(0,T;L^2(\Om))$ such that
$$\nabla_x P=\int_0^tf(s,x)ds-\rho^l\big(v-v^0)+2\mu\, {\rm div}_xe_x\Big(\int_0^tv(s,x)\,ds\Big).$$
Taking $p=\partial_t P\in W^{-1,\infty}(0,T;L^2(\om))$, we have that the second equation in (\ref{promo2}) holds and we formally write
$$P(t,x)=\int_0^tp(s,x)\,ds.$$
\end{Rem}
\section{ The viscoelastic model.}\label{secvimo}
 In the present section we prove that the inclusions of small liquid bubbles in an elastic medium provide a viscoelastic material.\par
We denote by $\Om$ a bounded open set in $\RR^N$ and by $O$ a connected open set with $\overline O\subset Y$, whose boundary is locally the graph of a Lipschitz function. The  coordinate system is chosen in such way that the null vector agrees with the center of mass of $\partial O$, i.e. such that
\beq\label{condO}\int_{\partial O}y\,d\sigma(y)=0.\eeq
Then, we define 
\beq\label{deOeHe} K_\ep=\big\{k\in\ZZ^N: \ C_k^\ep\subset \Om\big\},\quad \Om^l_\ep=\bigcup_{k\in K_\ep} (\ep k+\ep O),\quad \Om_\ep^s=\Om\setminus \overline{\Om^l_\ep},\quad \widetilde\Om_\ep=\bigcup_{k\in K_\ep}C_k^\ep.\eeq
\beq\label{deOmd} \Om^\delta=\big\{x\in\Om:\ {\rm dist}(x,	\partial\Om)>\delta\big\},\quad \forall\,\delta>0.\eeq
We introduce ${\cal P}:L^2(\partial O)^N\to L^2(\partial O)^N$ as the orthogonal projection from $L^2(\partial O)^N$ on the space of affine functions with a skew symmetric matrix (rigid movements), i.e. the space of functions $z=z(y)$ such that
$$z(y)=a+By,\quad \forall\, y\in \partial O,\quad a\in \RR^N,\ B\in \RR^{N\times N}\ \hbox{ skew-symmetric}.$$
We also take ${\cal Q}=I-{\cal P}$, the projection on the orthogonal space of the rigid movements.
\par\medskip
For a tensor $A\in {\cal L}(\RR^{N\times N}_s)$ which satisfies (\ref{proeA}), 
three positive constants $\rho^s$, $\rho^l$, $\mu$ and
$$f_\ep\in L^1(0,T;L^2(\Om))^N,\quad u_\ep^0\in H^1_{\partial\Om}(\Om^s_\ep)^N,\quad u_\ep^1\in L^2(\Om^s_\ep)^N,\quad v_\ep^0\in L^2(\Om_\ep^l)^N,$$
such that
\beq\label{hipcoi}
{\rm div}_x v_\ep^0=0\ \hbox{ in }\Om_\ep^l,\quad \sup_{\ep>0}\Big\{\big\| e(u_\ep^0)\big\|_{L^2(\Om_\ep^s)^{N\times N}}+\big\| u_\ep^1\big\|_{L^2(\Om_\ep^s)^N}+
\big\| v_\ep^0\big\|_{L^2(\Om_\ep^l)^N}\Big)<\infty,\eeq
\beq\label{hipfe} \exists f\in  L^1(0,T;L^2(\Om))^N,\ \hbox{ with }\ f_\ep\rightharpoonup f\ \hbox{ in }L^1(0,T;L^2(\Om))^N,\eeq
we are interested in the asymptotic behavior of the solutions of the linear system
\beq\label{sistho}\left\{\ba{l}\dis \rho^s\partial^2_{tt}u_\ep-{\rm div}_x(Ae_x(u_\ep))=f_\ep\ \hbox{ in }(0,T)\times\Om^s_\ep\\ \ecart\dis
\rho^l\partial_{t}v_\ep-2\mu\,{\rm div}_xe_x(v_\ep))+\nabla_x p_\ep=f_\ep\ \hbox{ in }(0,T)\times  \Om_\ep^l\\ \ecart\dis
{\rm div}_xv_\ep=0\ \hbox{ in }(0,T)\times  \Om_\ep^l\\ \ecart\dis \partial_t u_\ep=v_\ep,\ Ae_x(u_\ep)\nu=(2\mu\, e_x(v_\ep)-p_\ep I)\nu\quad\hbox{on }(0,T)\times \partial \Om_\ep^l\\ \ecart\dis
u_\ep=0\ \hbox{ on }(0,T)\times \partial\Om
\\ \ecart\dis (u_\ep)_{|t=0}=u_\ep^0,\quad (\partial_t u_\ep)_{|t=0}=u_\ep^1,\quad (v_\ep)_{|t=0}=v_\ep^0.  
\ea\right.
\eeq
By Theorem \ref{teoeu}, there exists a unique solution $(u_\ep,v_\ep,p_\ep)$ of (\ref{pbvsuv}), with
\beq\label{regsol} \ba{c}\dis u_\ep\in W^{1,\infty}(0,T;L^2(\Om_\ep^s))^N\cap L^\infty(0,T;H^1_{\partial\Om}(\Om_\ep^s))^N\\ \ecart\dis v_\ep\in 
L^\infty(0,T;L^2(\Om_\ep^l))^N\cap L^2(0,T;H^1(\Om_\ep^l))^N\\ \ecart\dis
p_\ep\in W^{-1,\infty}(0,T; L^2(\Om_\ep^l)).\ea\eeq
Moreover, it satisfies the energy identity
\beq\label{idente}\ba{l}\dis
{1\over 2}{d\over dt}\left(\rho^s\int_{\Om_\ep^s}|\partial_t u_\ep|^2dx+\int_{\Om_\ep^s}Ae_x(u_\ep):e_x(u_\ep)dx+\rho^l\int_{\Om^l_\ep}|v_\ep|^2dx\right)\\ \ecart\dis\qquad \qquad+2\mu\int_{\Om^l_\ep}|e_x( v_\ep)|^2dx=\int_{\Om_\ep^s}f_\ep\cdot \partial_t u_\ep\,dx+\int_{\Om_\ep^l}f_\ep\cdot v_\ep\,dx\quad \hbox{ in }(0,T).\ea\eeq
Let us also use the following results:\par\medskip
Next Lemma can be found for example in \cite{OlShYo}.
\begin{Lem} \label{opext} There exists a sequence of continuous linear extension operators ${\cal E}_\ep:H^1_{\partial\Om}(\Om_\ep^s)^N\to H^1_0(\Om)^N$, whose norm is bounded independently of $\ep$, i.e. they satisfy 
\beq\label{EepoP}{\cal E}_\ep u=u\ \hbox{ in }\Om_\ep^s,\quad
\big\|{\cal E}_\ep u\big\|_{H^1_0(\Om)^N}\leq C\|e_x(u)\|_{L^2(\Om_\ep^s)^{N\times N}},\quad \forall\, u\in H^1_{\partial\Om}(\Om_\ep^s)^N,\eeq
with $C>0$ independent of $u$ and $\ep$.
\end{Lem} 
\begin{Lem}  \label{leext2}
There exists $C>0$ such that for every $k\in\ZZ^N$, every $\ep>0$, and every $p\in L^2(\ep (k+O))$, there exits $w\in H^1_0(C^k_\ep)$ satisfying
\beq\label{leacpr}{\rm div}_x w=p\ \hbox{ in }\ep(k+O),\quad \|w\|_{H^1_0(C^k_\ep)^N}\leq C\|p\|_{L^2(\ep(k+O))}.\eeq
\end{Lem}\par\medskip\noindent
{\bf Proof}. By the regularity of $O$,  there exists $C>0$, such that for every $\hat p\in L^2(O)$, there exists $\hat w\in H^1_0(Y)^N$ satisfying
$${\rm div}_y \hat w=\hat p\ \hbox{ in }O,\quad  \|w\|_{H^1_0(Y)^N}\leq C\|\hat p\|_{L^2(O)}.$$
The proof of (\ref{leacpr}) then follows by using the change of variables $y=(x-\ep k)/\ep$ which transforms $\ep(k+O)$ in $O$ and
$C_k^\ep$ in $Y$. \cqfd
\par\medskip
As a consequence of (\ref{idente}) and Lemmas \ref{opext} and \ref{leext2}, we get
\begin{Pro} There exists $C>0$, which does not depend on $\ep$, the right-hand side $f_\ep$ and the initial conditions $u^0_\ep$, $u^1_\ep$, $v_\ep$, such that the solution $(u_\ep,v_\ep,p_\ep)$ of (\ref{sistho}) satisfies
\beq\label{estpriso} \ba{l}\dis\|u_\ep\|_{W^{1,\infty}(0,T;L^2(\Om_\ep^s))^N\cap L^\infty(0,T;H^1_{\partial\Om}(\Om_\ep^s))^N}\\ \ecart\dis+\|v_\ep\|_{L^\infty(0,T;L^2(\Om_\ep^l))^N}
+\|e_x(v_\ep)\|_{L^2(0,T;L^2(\Om_\ep^l))^{N\times N}}
+\Big\|\int_0^tp_\ep(s,x)\,ds\Big\|_{L^\infty(0,T;L^2(\Om_\ep^l))}\\ \ecart\dis\leq C\Big(\|e_x(u_\ep^0)\|_{L^2(\Om^s_\ep)^{N\times N}}+\|u_\ep^1\|_{L^2(\Om^s_\ep)^N}+\|v^0_\ep\|_{L^2(\Om^l_\ep)^N}+\|f_\ep\|_{L^1(0,T;L^2(\Om))^N}\Big).\ea\eeq
\end{Pro}
\par\medskip\noindent
{\bf Proof.} By (\ref{idente}) and Gronwall's inequality, we have
$$\ba{l}\dis\|u_\ep\|_{W^{1,\infty}(0,T;L^2(\Om_\ep^s))^N}+\|e_x(u_\ep)\|_{L^\infty(0,T;L^2(\Om_\ep^s))^{N\times N}}+\|v_\ep\|_{L^\infty(0,T;L^2(\Om_\ep^l))^N}
+\|e_x(v_\ep)\|_{L^2(0,T;L^2(\Om_\ep^l))^{N\times N}}
\\ \ecart\dis\leq C\Big(\|e_x(u_\ep^0)\|_{L^2(\Om^s_\ep)^{N\times N}}+\|u_\ep^1\|_{L^2(\Om^s_\ep)^N}+\|v^0_\ep\|_{L^2(\Om^l_\ep)^N}+\|f_\ep\|_{L^1(0,T;L^2(\Om))^N}\Big),\ea$$
where by Lemma \ref{opext} and Korn's inequality, we also have
$$\|u_\ep\|_{H^1_{\partial\Om}(\Om_\ep^s)^N}\leq \|{\cal E}_\ep u_\ep\|_{H^1_{\partial\Om}(\Om)^N}\leq C\|e_x({\cal E}u_\ep)\|_{L^2(\Om)^{N\times N}}\leq  C\|e_x(u_\ep)\|_{L^2(\Om_\ep^s)^{N\times N}}.$$
The estimate for the integral of $p_\ep$ follows from Lemma \ref{leext2} reasoning as in the proof of (\ref{estpri2}).\par\hfill \cqfd\par\medskip
The following proposition provides a limit for the initial conditions. Observe that thanks to Lemma \ref{opext} we can assume $u_\ep^0$ defined in the whole of $\Om$.
\begin{Pro} \label{Prolic} There exists a subsequence of $\ep$, still denoted by $\ep$,   $u^0\in H^1_0(\Om)^N$ and $u^1,v^0\in L^2(\Om)^N$ such that 
\beq\label{convduep0} u_\ep^0\rightharpoonup u^0\ \hbox{ in }H^1_0(\Om)^N,\qquad u_\ep^1\chi_{\Om_\ep^s}\rightharpoonup (1-|O|)u^1,\ v^\ep\chi_{\Om_\ep^l}\rightharpoonup |O|v^0\quad\hbox{ in }L^2(\Om)^N.\eeq
Moreover, defining $\hat u_\ep^0$ in $\Om\times (Y\setminus\overline O)$  by
\beq\label{dedai2e} \hat u_\ep^0(x,y)=u_\ep^0\Big(\ep\kappa\big({x\over\ep}\big)+\ep y\Big),\ \hbox{ a.e. }(x,y)\in \Om\times (Y\setminus\overline O),
\eeq
there exist $u^0_1\in L^2(\Om;H^1_\sharp(Y\setminus\overline O))^N$, such that
\beq\label{limcoiu0} {1\over\ep}D_y\hat u^0_\ep\rightharpoonup D_xu^0(x) y+D_yu^0_1\ \hbox{ in }L^2(\Om;L^2(Y\setminus \overline O))^{N\times N}.\eeq
\end{Pro}\par\medskip\noindent
{\bf Proof.} Assertion (\ref{convduep0}) just follows from $u_\ep^0$ bounded in $H^1_0(\Om)^N$ and $u_\ep^1\chi_{\Om_\ep^s}$, $v^\ep\chi_{\Om_\ep^l}$ bounded in $L^2(\Om)^N$.\par
Assertion (\ref{limcoiu0}) is also a classical result for the unfolding method (\cite{CiDaGr}), which is strongly related to the classical compactness results in the two-scale convergence theory for sequences bounded in $H^1(\Om)$ (\cite{All}, \cite{Ngu}). Let us briefly recall the idea of the proof because we will need it to prove Theorem \ref{Thhocde} below. \par
Taking into account that
$${1\over \ep^2}\into\int_{Y\setminus \overline O}|D_y\hat u^0_\ep|^2dx=\int_{\Om_\ep^s}|D_x u^0_\ep|^2dx\leq C,$$
we introduce
\beq\label{otfu}\breve u_\ep^0={1\over\ep}\Big(\hat u_\ep^0-\intm_{\partial O}\hat u_\ep^0\,d\sigma(\beta)\Big).
\eeq
Thanks to Poincar\'e-Wirtinger's inequality, $\check u_\ep^0$ is a bounded sequence in $L^2(\Om;H^1_\sharp(Y\setminus\overline O))^N$. Thus, for another subsequence, there exists $\breve u^0\in L^2(\Om;H^1_\sharp(Y\setminus\overline O))^N$ such that
\beq\label{conbu}\breve u_\ep^0\rightharpoonup \breve u^0\ \hbox{ in }L^2(\Om;H^1_\sharp(Y\setminus\overline O))^N.\eeq
Now, we observe that by definition of $\breve u_\ep^0$, for a.e. $x\in \Om$ and a.e. $y'\in (-1/2,1/2)^{N-1}$, we have
$$\hat u_\ep^0\Big(x+\ep e_1,-{1\over 2},y'\Big)=\hat u_\ep^0\Big(x,{1\over 2},y'\Big).$$
Therefore, the sequence $\check u_\ep^0$ satisfies
\beq\label{igprpe}\breve u_\ep^0\Big(x+\ep e_1,-{1\over 2},y'\Big)-\breve u_\ep^0\Big(x,{1\over 2},y'\Big)=
-\intm_{\ep\kappa({x\over\ep})+\ep\partial O}\hskip-10pt {u^0_\ep(\beta+\ep e_1)-u^0_\ep(\beta)\over\ep}d\sigma(\beta).\eeq
In this equality, the left-hand side converges weakly in $L^2(\Om;H^{1\over 2}((-1/2,1/2)^{N-1}))^N$ to
$$\breve u^0\Big(x,-{1\over 2},y'\Big)-\breve u^0\Big(x,{1\over 2},y'\Big),$$
while for the right-hand side, using that $u_\ep^0$ it is bounded in $H^1_0(\Om)^N$, it is simple to check that
$$\intm_{\ep\kappa({x\over\ep})+\ep\partial O}\hskip-10pt {u^0_\ep(\beta+\ep e_1)-u^0_\ep(\beta)\over\ep}d\sigma(z)\rightharpoonup \partial_1u^0\ \hbox{ in }{\cal D}'(\Om).$$
Thus, we get
$$\breve u^0\Big(x,-{1\over 2},y'\Big)-\breve u^0\Big(x,{1\over 2},y'\Big)=-\partial_1u^0(x)$$
Reasonig analogously with the other components we deduce that
\beq\label{deu10}u_1^0(x,y):=\breve u^0(x,y)-D_x u^0(x) y\eeq
is periodic with respecto to $y$ of period $Y$ and by (\ref{otfu}) and (\ref{conbu}), it satisfies (\ref{limcoiu0}). \cqfd\par\medskip
Theorem \ref{Thhocde} below provides a partial differential problem for the limits of the solutions of (\ref{sistho}) where both the macroscopic variable $x$ and the microscopic variable $y$ intervene (two-scale problem). Taking into account Lemma \ref{opext} and extending $u_\ep$ by zero outside of $\Om$, we assume that $u_\ep$ is defined in
$(0,T)\times\RR^N$. By Proposition \ref{Prolic}, it is not restrictive to assume that the initial conditions satisfy (\ref{convduep0}) and (\ref{limcoiu0}).
\begin{Thm}\label{Thhocde} We assume that the functions $u_\ep^0,u_\ep^1$ and $v_\ep^0$ in (\ref{sistho}) satisfy (\ref{convduep0}) and (\ref{limcoiu0}) for some $u^0\in H^1_0(\Om)^N$, $u^1,v^0\in L^2(\Om)^N$ and $u_1^0\in L^2(\Om;H^1_\sharp(Y\setminus\overline O))^N$. \par
For $(u_\ep,v_\ep,p_\ep)$ the solution of (\ref{sistho}), we define $\hat u_\ep$ in  $(0,T)\times\Om\times (Y\setminus \overline O)$, $\hat v_\ep,\hat p_\ep$ in  $(0,T)\times\tilde \Om_\ep\times O$ 
by
\beq\label{defhuhv} \ba{c}\dis\hat u_\ep(t,x,y)=u_\ep\Big(t,\ep\kappa\big({x\over\ep}\big)+\ep y\Big),\quad \hat v_\ep(t,x,y)=v_\ep\Big(t,\ep\kappa\big({x\over\ep}\big)+\ep y\Big)\\ \ecart\dis \hat p_\ep(t,x,y)=p_\ep\Big(t,\ep\kappa\big({x\over\ep}\big)+\ep y\Big).\ea\eeq
Defining
\beq\label{defrefe} \rho_{eff}:=\rho^s(1-|O|)+\rho^l|O|,\eeq
we also take
\beq\label{terpbli}\ba{c} u\in W^{1,\infty}(0,T;L^2(\Om))^N\cap L^\infty(0,T;H^1_0(\Om))^N,\quad
u_1\in L^\infty(0,T;L^2(\Om;H^1_\sharp(Y\setminus\overline O)))^N\\ \ecart\dis
v\in L^\infty(0,T;L^2(\Om;H^1(O)))^N,\quad p\in L^\infty(0,T;L^2(\Om;L^2(O))),\ea\eeq
a solution of 
\beq\label{pbli2e} \left\{\ba{l}\dis \rho_{eff}\partial^2_{tt}u-{\rm div}_x\Bigg(\int_{Y\setminus\overline O}A\big(e_x(u)+e_y(u_1)\big)dy\\ \ecart\dis\qquad\qquad\qquad\qquad+\int_O(2\mu e_y(v)-pI)dy\Bigg)=f\ \hbox{ in }(0,T)\times\Om\\ \ecart\dis
-{\rm div}_y\big(Ae_y(u_1)\big)=0\ \hbox{ in }(0,T)\times \Om\times \big(Y\setminus \overline O)\\ \ecart\dis
-2\mu\,{\rm div}_y e_y(v)+\nabla_yp=0\ \hbox{ in }(0,T)\times\Om\times O\\ \ecart\dis
{\rm div}_y v=0\ \hbox{ in } (0,T)\times \Om\times O\\ \ecart\dis
u=0\ \hbox{ on }(0,T)\times\partial\Om\\ \ecart\dis
u_1\hbox{ periodic in }Y,\quad A\big(e_x(u)+e_y(u_1)\big)\nu\ \hbox{anti-periodic in }Y,\\ \ecart\dis
{\cal Q}\big(e_x(\partial_tu)y+\partial_tu_1\big)=v\ \hbox{ on }(0,T)\times\Om\times\partial O\\ \ecart\dis
{\cal P}\big(A\big(e_x(u)+e_y(u_1)\big)\nu\big)=0\ \hbox{ on }(0,T)\times\Om\times\partial O\\ \ecart\dis
A\big(e_x(u)+e_y(u_1)\big)\nu=(2\mu e_y(v)-pI)\nu\ \hbox{ on }(0,T)\times \Om\times\partial O\\ \ecart\dis
u_{|t=0}=u^0,\  \rho_{eff}\partial_tu_{|t=0}=\rho^s(1-|O|)u^1+\rho^l|O| v^0\quad\hbox{ in }\Om\\ \ecart\dis
({\cal Q}u_1)_{|t=0}={\cal Q}u_1^0\ \hbox{ on }\Om\times\partial O.
\ea\right.\eeq
Then, we have
\beq\label{convso} \left\{\ba{l}\dis u_\ep\stackrel{\ast}\rightharpoonup u\ \hbox{ in }W^{1,\infty}(0,T; L^2(\Om))^N\cap L^\infty(0,T;H^1_0(\Om))^N\\ \ecart\dis
v_\ep\chi_{\Om_\ep^l}\stackrel{\ast}\rightharpoonup |O|\partial_t u\ \hbox{ in }L^\infty(0,T;L^2(\Om))^N\\ \ecart\dis {1\over\ep}D_y\hat u_\ep\stackrel{\ast}\rightharpoonup D_xu+D_yu_1\ \hbox{ in }L^\infty(0,T;L^2(\Om;L^2(Y\setminus O)))^{N\times N}\\ \ecart\dis
{1\over \ep}e_y(\hat v_\ep)\rightharpoonup e_y(v)\ \hbox{ in }L^2(0,T;L^2(\Om^\delta;H^1(O)))^{N\times N},\ \forall\, \delta>0\\ \ecart\dis
\int_0^t\hat p_\ep\,dt\stackrel{\ast}\rightharpoonup \int_0^tp\,dt\ \hbox{ in }L^\infty(0,T;L^2(\Om^\delta;L^2(O))),\ \forall\, \delta>0.
\ea\right.
\eeq
Moreover, the functions $u$, $u_1$ and $v$ satisfy the energy identity
\beq\label{idenpl}\ba{l}\dis {1\over 2}{d\over dt}\left(\into |\partial_tu|^2dx+\into\int_{Y\setminus\overline O} A\big(e_x(u)+e_y(u_1)\big):\big(e_x(u)+e_y(u_1)\big)dydx\right)\\ \ecart\dis
+2\mu\into\int_O|e_y(v)|^2dydx=\into f\cdot\partial_tu\,dx\ \hbox{ in }(0,T).\ea\eeq
\end{Thm}
\par\medskip\noindent
{\bf Proof.} We divide it in three steps.\par\medskip\noindent
{\it Step 1.} Since $u_\ep$ is bounded in $W^{1,\infty}(0,T; L^2(\Om))^N\cap L^\infty(0,T;H^1_0(\Om))^N$, there exist a subsequence of $\ep$, still denoted by $\ep$, and $u\in W^{1,\infty}(0,T; L^2(\Om))^N\cap L^\infty(0,T;H^1_0(\Om))^N$ such that the first line in (\ref{convso}) holds.
From (\ref{sistho}) and (\ref{convduep0}), the function $u$ also satisfies
\beq\label{coninu} u_{|t=0}=u^0\ \hbox{ in }\Om.\eeq\par
Now, using the change of variables $y=(x-\ep k)/\ep$ in each cube $C_k^\ep$, we deduce from estimate (\ref{estpriso}) and the definitions of $\hat u_\ep$, $\hat v_\ep$ and $\hat p_\ep$ that
\beq\label{acohu}\|\hat u_\ep\|_{W^{1,\infty}(0,T;L^2(\RR^N; L^2(Y\setminus \overline O)))^N}+{1\over\ep}\|D_y\hat u_\ep\|_{L^\infty(0,T;L^2(\RR^N; L^2(Y\setminus O)))^{N\times N}}\leq C,\eeq
\beq\label{acohv}\|\hat v_\ep\|_{L^{\infty}(0,T;L^2(\tilde \Om_\ep; L^2( O)))^N}+{1\over\ep}\|e_y(\hat v_\ep)\|_{L^2(0,T;L^2(\tilde\Om_\ep; L^2(O)))^{N\times N}}\leq C,\eeq
\beq\label{acohp}\left\|\int_0^t\hat p_\ep\,ds\right\|_{L^\infty(0,T;L^2(\tilde \Om_\ep;L^2(O)))}\leq C.\eeq
We introduce 
\beq\label{decucv}\check u_\ep={1\over\ep}\Big(\hat u_\ep-\intm_{\partial O} \hat u_\ep d\sigma(z)\Big)\eeq\par
From (\ref{acohu}) and Korn's inequality,  $\check u_\ep$ is bounded in $L^\infty(0,T;L^2(\Om;H^1(Y\setminus O))^N$. Reasoning as in the proof of Proposition \ref{Prolic} we deduce from this estimate the existence of  $u_1\in L^\infty(0,T;L^2(\Om;H^1_\sharp(Y\setminus \overline O)))^N$ such that, for a subsequence, the third convergence in (\ref{convso}) holds.
\par Analogously, we take
\beq\label{dechv}\check v_\ep= {1\over\ep}{\cal Q}\hat v_\ep.\eeq
By Korn's inequality  $\|\check v_\ep\|_{L^2(0,T;L^2(\tilde\Om_\ep;H^1(O)))^N}$ is bounded. Thus, for a subsequence, there exists $v\in L^2(0,T;L^2(\Om;H^1(O)))^N$ such that
\beq\label{concve} \check v_\ep\rightharpoonup v\ \hbox{ in }L^2(0,T;L^2(\Om^\delta;H^1(O)))^N,\quad\forall\,\delta>0.\eeq
In particular this proves the fourth convergence in (\ref{convso}).\par
Taking into account that $\partial_t u_\ep=v_\ep$ on $(0,T)\times \partial\Om_\ep^l$, and the definitions of $\check u_\ep$ and $\check v_\ep$ we get
\beq\label{igufoh}{\cal Q}\partial_t\check u_\ep=\check v_\ep\ \hbox{ on }(0,T)\times\tilde \Om_\ep\times\partial O,\eeq
and then, by  (\ref{convso}), line 8 in (\ref{pbli2e}). By (\ref{igufoh}), we also have that the sequence $\check u_\ep$ is bounded in $H^1(0,T;L^2(\Om;H^{1\over 2}(\partial O)))^N$. Using then that $(\check u_\ep)_{|t=0}=\check u_\ep^0,$ with $\check u_\ep^0,$ defined by (\ref{otfu}),
and (\ref{conbu}), (\ref{deu10}) we deduce the last assertion in (\ref{pbli2e}).\par
Since ${\rm div}_x v_\ep=0$ in $(0,T)\times \Om$ implies ${\rm div}_y\hat v_\ep=0$ in $(0,T)\times \tilde\Om_\ep\times O$ and  the rigid movements have null divergence, we get from definition (\ref{dechv}) of $\check v_\ep$ and (\ref{concve}) that ${\rm div}_yv=0$ in $(0,T)\times\Om\times O$.\par
Now, we observe that $\|\check v_\ep\|_{L^2(0,T;L^2(\tilde\Om_\ep;H^1(O)))^N}$ bounded and definition (\ref{dechv}) of $\hat v_\ep$ imply 
$${\cal Q}\hat v_\ep=\hat v_\ep-{\cal P}\hat v_\ep\rightarrow 0\ \hbox{ in }L^2(0,T;L^2(\Om^\delta;H^1(O)))^N,\quad\forall\,\delta>0.$$
On the other hand,  the first convergence in (\ref{convso}) implies that $u_\ep$ converges strongly to $u$ in $L^\infty(0,T;L^2(\Om))^N$ and then, by (\ref{acohu}) that
$$\hat u_\ep\stackrel{\ast}\rightharpoonup u\ \hbox{ in }W^{1,\infty}(0,T; L^2(\Om;L^2(Y\setminus \overline O)))^N\cap L^\infty(0,T; L^2(\Om;H^1(Y\setminus \overline O)))^N.$$
Using that $\partial_t\hat u_\ep=\hat v_\ep$ on $(0,T)\times \tilde\Om_\ep\times \partial O$, 
 we conclude that 
$$\hat v_\ep=\big(\hat v_\ep-{\cal P}v_\ep\big)+{\cal P}v_\ep
\rightharpoonup {\cal P}\partial_t u=\partial_t u\ \hbox{ in }L^2(0,T;L^2(\Om^\delta;H^1(O)))^N,$$
where we have used that ${\cal P}\partial_t u=\partial_t u$ because $\partial_tu$ does not depend on the microscopic variable $y$. This proves the second convergence in (\ref{convso}).\par
By (\ref{acohp}), we also deduce the existence of $P\in L^\infty(0,T;L^2(\Om;L^2(O)))$, such that
\beq\label{convppe}\int_0^t\hat p_\ep\,dt\stackrel{\ast}\rightharpoonup P\ \hbox{ in }L^\infty(0,T;L^2(\Om^\delta;L^2(O))),\ \forall\, \delta>0.\eeq
\par\medskip\noindent
{\it Step 2.}  
For
\beq\label{hipfute}\left\{\ba{c}\dis \varphi\in C^\infty([0,T]; {\cal D}(\Om))^N,\quad \phi\in C^\infty([0,T]; {\cal D}(\Om;C^\infty_\sharp(\overline{Y\setminus O})))^N\\ \ecart\dis
\psi\in C^\infty([0,T];{\cal D}(\Om;C^\infty_\sharp(\overline O)))^N\\ \ecart\dis
{\cal Q}(e_x(\varphi)y+\phi)=\psi\ \hbox{ on }(0,T)\times \Om\times \partial O,\quad \varphi, \phi,\psi\ \hbox{ vanishing for }t=T,\ea\right.\eeq
we take 
$$\varphi(t,x)+\ep\phi\big(t,x,{x\over\ep}\big)$$
as test  function in the first equation in  (\ref{sistho}), and
$$\varphi(t,x)+\ep\psi\Big(t,x,{x\over\ep}\Big)+\ep{\cal P}\phi\big(t,x,{x\over\ep})
-{\cal Q}\Big(e_x(\varphi)\big(x-\ep\kappa\big({x\over\ep}\big)\big)\Big)$$
as test function in the second one. Adding the corresponding equalities, using that these test functions agree on $\partial\Om_\ep^l$ and that the rigid movements have null divergence, we get
$$\ba{l}\dis
-\rho^s\int_{\Om_\ep^s}u_\ep^1\cdot\varphi(0,x)dx
-\rho^l\int_{\Om_\ep^l}v^0_\ep\cdot\varphi(0,x)dx
-\rho^s\int_0^T\hskip-5pt\int_{\Om_\ep^s}\partial_t u_\ep\cdot \partial_t\varphi dxdt-\rho^l\int_0^T\hskip-5pt\int_{\Om_\ep^l}v_\ep\cdot \partial_t\varphi dxdt\\ \ecart\dis
+\int_0^T\hskip-5pt\int_{\Om_\ep^s}Ae_x(u_\ep):\Big(e_x(\varphi)+e_y(\phi)\big(t,x,{x\over\ep}\big)\Big)dxdt+2\mu\int_0^T\hskip-5pt\int_{\Om_\ep^l}e_x(v_\ep): e_y(\psi)\Big(t,x,{x\over\ep}\Big)dxdt\\ \ecart\dis+\int_0^T\hskip-5pt\int_{\Om_\ep^l} \int_0^tp_\ep(s,x,y)ds\, {\rm div}_y\partial_t\psi\big(t,x,{x\over\ep}\big)dxdt=\int_0^T\into f\cdot\varphi\,dxdt+O_\ep,\ea$$ \par
In the three last terms, we use the change of variables $y=(x-\ep k)/\ep$ in every cube $C_k^\ep$. Then, by  (\ref{convduep0}) and (\ref{convso}), we can pass to the limit in this equality to deduce
\beq\label{pblim}\ba{l}\dis
-\into\left(\rho^su^1(1-|O|)+\rho^lv^0|O|\right)\cdot\varphi(0,x)dx
-\rho_{eff}\int_0^T\hskip-5pt\into\partial_t u\cdot \partial_t\varphi\,dxdt\\ \ecart\dis
+\int_0^T\hskip-5pt\into\int_{Y\setminus O}A\big(e_x(u)+e_y(u_1)\big):\big(e_x(\varphi)+e_y(\phi)\big)dydxdt\\ \ecart\dis
+2\mu\int_0^T\hskip-5pt\into\int_Oe_y(v):e_y(\psi)dxdt+\int_0^T\hskip-5pt\into\int_O P\, {\rm div}_y\partial_t\psi\,dydxdt\\ \ecart\dis=\int_0^T\hskip-5pt\into f\cdot\varphi\,dxdt.\ea\eeq
By density this equality holds true just assuming
\beq\label{hipfuteD}\left\{\ba{c}\dis \varphi\in W^{1,1}(0,T; L^2(\Om))^N\cap L^1(0,T;H^1_0(\Om))^N,\quad \phi\in L^1(0,T;L^2(\Om;H^1_\sharp(Y\setminus \overline O)))^N\\ \ecart\dis
\psi\in L^1(0,T;L^2(\Om;L^2(O)))^N,\quad {\rm div}_y\partial_t\psi\in L^1(0,T; L^2(\Om;L^2(O)))\\ \ecart\dis
{\cal Q}(e_x(\varphi)y+\phi)=\psi\ \hbox{ on }(0,T)\times \Om\times \partial O\\ \ecart\dis
 \varphi_{|t=T}=0\ \hbox{ in }\Om,\quad ({\rm div}_y\psi)_{|t=T}=0\ \hbox{ in }\Omega\times O.\ea\right.\eeq
This is equivalent to (\ref{pbli2e}), with $p$ and $P$ related by (see Remark \ref{reep})
$$\int_0^t p(s,x,y)ds=P(t,x,y)\ \hbox{ a.e. }(t,x,y)\in (0,T)\times \Om\times O,$$
where the first equation is obtained taking in (\ref{pblim}) $\phi=0$, $\psi={\cal Q}(e_x(\varphi)y),$ and using that  lines 3, 9 and 10 in (\ref{pbli2e}) imply
$$\ba{l}\dis \int_0^T\int_O \big(2\mu e_y(v)-pI\big):e_y\big({\cal Q}(e_x(\varphi)y)\big)dydt=
\int_0^T\int_O \big(2\mu e_y(v)-pI\big)dy:e_x(\varphi)dt\\ \ecart\dis=
2\mu\int_0^T\int_O  e_y(v)dy:e_x(\varphi)dt+\int_0^T\int_O\int_0^tp(s,x,y)dsdy:{\rm div}_x\partial_t\varphi\,dt.\ea$$
The regularity of $p$ from lines 4 and 10 in (\ref{pbli2e}).
 \par\medskip\noindent
{\it Step 3.}
If we assume $\partial_tu\in L^1(0,T; H^1_0(\Om))^N$ and $\partial_tu_1\in L^1(0,T; L^2(\Om;H^1_\sharp(Y\setminus\overline O)))$, we can multiply the first equation in (\ref{pbli2e}) by $\partial_tu$ and to integrate in $\Om$, the second one by $\partial_tu_1$ and to integrate in $\Om\times (Y\setminus\overline O)$, and the third one by $v$ and to integrate in $\Om\times O$. Adding the corresponding inequalities
we deduce (\ref{idenpl}). The proof of the result when $u$, $u_1$ are just in $W^{1,1}(0,T;L^2(\Om))^N\cap L^\infty(0,T; H^1_0(\Om))^N$ and $L^\infty(0,T; L^2(\Om;H^1_\sharp(Y\setminus\overline O)))^N$ respectively follows by a regularization argument.\par
The uniqueness result for (\ref{pbli2e}) mentioned in Remark \ref{remuni} below implies that there is no
 need to extract any subsequence in (\ref{convso}).
 \cqfd
 \begin{Rem} \label{remuni} From (\ref{idenpl}) and lines 4 and 10  in (\ref{pbli2e})
we get that $u$, $v$, $p$ are uniquely defined while $u_1$ is  defined except for an adding term  $R(t,x)y+r(t,x)$, with $R$ skew-symmetric. In particular $e_x(u_1)$ is uniquely defined.
\end{Rem}
\begin{Rem} \label{remIEl} Related to (\ref{idenpl}) we observe that taking into account lines 3, 7, 9 and 12  in (\ref{pbli2e}), we get  
$$\ba{l}\dis\left(\into\int_{Y\setminus\overline O} A\big(e_x(u)+e_y(u_1)\big):\big(e_x(u)+e_y(u_1)\big)dydx\right)_{|t=0}\\ \ecart\dis =\into\int_{Y\setminus\overline O} A\big(e_x(u^0)+e_y(u^\ast_1)\big):\big(e_x(u^0)+e_y(u^\ast_1)\big)dydx,\ea$$
where  $u^\ast_1$ is the unique solution of
\beq\label{deu1a} \left\{\ba{l}\dis 
-{\rm div}_yA\big(e_x(u^0)+e_y(u_1^\ast)\big)=0\ \hbox{ in }\Om\times(Y\setminus \overline O)
\\ \ecart\dis
u^\ast_1\hbox{ periodic in }Y,\ A\big(e_x(u^0)+e_y(u^\ast_1)\big)\nu \hbox{ anti-periodic in }Y\\ \ecart\dis
{\cal Q}u^\ast_1={\cal Q}u_1^0\ \hbox{ on }\Om\times \partial O\\ \ecart\dis
{\cal P}\big(A\big(e_x(u^0)+e_y(u^\ast_1)\big)\nu\big)=0\ \hbox{ on }\Om\times\partial O.
\ea\right.\eeq
\end{Rem}

 \begin{Rem} Problem (\ref{pbli2e}) can also be written as the variational system given by (\ref{terpbli}), (\ref{pblim}), (\ref{hipfuteD}) and
\beq\label{capvl}\ba{c}\dis  {\rm div}_yv=0\ \hbox{ in }(0,T)\times \Om\times O\\ \ecart\dis
{\cal Q}\big(e_x(\partial_tu)y+\partial_tu_1\big)=v\hbox{ on }(0,T)\times \Om\times\partial O.\ea\eeq
\end{Rem}
\par\medskip
Our purpose now is to obtain an equation for the function $u$ in (\ref{pbli2e}), which will describe the macroscopic behavior of the solutions of (\ref{sistho}). For this aim, we need to eliminate the functions $u_1$, $v$ and $p$ in (\ref{pbli2e}), depending on the microscopic variable $y$. Observe that for a.e. $x\in \Om$, they can be obtained from $e_x(u)(x,.)$, solving the system corresponding to eliminate  lines  1, 2, 6 and 11 in (\ref{pbli2e}). Let us study the existence and uniqueness of solution for this problem in the following two lemmas. They will also provide  some qualitative properties for these functions.
\begin{Lem}\label{leexspme} For every $\zeta^0\in H^{1\over 2}(\partial O)^N$, with ${\cal P}\zeta^0=0$, there exits a unique solution $(\zeta,\eta,q)$ of 
\beq\label{pbvame}\left\{\ba{l}\dis \zeta\in C^\infty([0,\infty);H^1_\sharp (Y\setminus\overline O))^N,\ \ 
\eta\in C^\infty([0,\infty);H^1(O))^N,\ \ q\in C^\infty([0,\infty);L^2(O))\\ \ecart\dis 
-{\rm div}_y\big(Ae_y(\zeta)\big)=0\ \hbox{ in }(0,\infty)\times(Y\setminus\overline O)\\ \ecart\dis
-2\mu\,{\rm div}_ye_y(\eta)+\nabla_y q=0\ \hbox{ in }(0,\infty)\times O
\\ \ecart\dis {\rm div}_y\eta=0\ \hbox{ in }(0,\infty)\times O\\ \ecart\dis
\zeta\hbox{ periodic in }Y,\ \ Ae_y(\zeta)\nu\hbox{ anti-periodic in }Y
\\ \ecart\dis 
{\cal Q}\partial_t\zeta=\eta,\  {\cal P}\big(Ae_y(\zeta)\nu\big)=0,\ Ae_y(\zeta)\nu=\big(2\mu e_y(\eta)-q I\big)\nu\ \hbox{ on }(0,\infty)\times \partial O\\ \ecart\dis
({\cal Q}\zeta)_{|t=0}=\zeta^0\ \hbox{ on }\partial O.
\ea\right.
\eeq
Moreover, it satisfies the energy identity
\beq\label{ideeL1}{1\over 2}{d\over dt}\int_{Y\setminus \overline O}Ae_y(\zeta):e_y(\zeta)dy+
2\mu\int_O|e_y(\eta)|^2dy=0\ \hbox{ in }[0,\infty),\eeq
and there exist two constants $\lambda,\Lambda>0$ such that
\beq\label{decexpl}\|\zeta\|_{H^1_\sharp(Y\setminus\overline O)^N}+\|\eta\|_{H^1(O)^N}+\|q\|_{L^2(O)}\leq \Lambda e^{-\lambda t}\|\zeta^0\|_{H^{1\over 2}(\partial O)^N}\ \hbox{ in }(0,\infty).\eeq
\end{Lem}
\par\noindent
{\bf Proof.} We introduce the Hilbert space
$${\cal H}=\big\{h\in H^{1\over 2}(\partial O)^N:\ {\cal P}h=0\big\}.$$
Then, we define ${\cal R}=({\cal R}_1,{\cal R}_2,{\cal R}_3): {\cal H}\to H^1_\sharp(Y\setminus \overline O)^N\times H^1(O)^N\times L^2(O)$
in the following way: \par
For $h\in H$, ${\cal R}_1h$ is the solution of the linear elasticity problem
\beq\label{defR1} \left\{\ba{l}\dis -{\rm div}_y\big(Ae_y({\cal R}_1h)\big)=0\ \hbox{ in }Y\setminus\overline O\\ \ecart\dis
{\cal R}_1h\hbox{ perodic in Y},\ Ae_y({\cal R}_1h)\nu\hbox{ anti-periodic in }Y\\ \ecart\dis
{\cal Q}{\cal R}_1h=h,\ {\cal P}\big(Ae_y({\cal R}_1h)\nu\big)=0\ \hbox{ on }\partial O,\quad \int_{\partial O}{\cal R}_1h\,d\sigma(x)=0.\ea\right.\eeq
The existence and uniqueness of this problem follows from Lax-Milgram's theorem  writing it in the variational form
$$\left\{\ba{l}\dis{\cal R}_1h\in H^1_\sharp(Y\setminus \overline O),\quad {\cal Q}{\cal R}_1h=h\ \hbox{ on }\partial O,\quad \int_{\partial O}{\cal R}_1h\,d\sigma(x)=0\\ \ecart\dis
\int_{Y\setminus \overline O}Ae_y({\cal R}_1h):e_y(\phi)dy=0,\quad
\dis \forall\,\phi\in H^1_\sharp (Y\setminus\overline O)^N,\  {\cal Q}\phi=0\ \hbox{ on }\partial O.\ea\right.$$
The pair $\big({\cal R}_2h,{\cal R}_3h)$ is defined as the unique  solution of the Stokes problem
\beq\label{defR2} \left\{\ba{l}\dis -2\mu\,{\rm div}_ye_y\big({\cal R}_2h\big)+\nabla_y {\cal R}_3h=0\ \hbox{ in }O\\ \ecart\dis
{\rm div}_y{\cal R}_2h=0\ \hbox{ in }O\\ \ecart\dis
\big(2\mu e_y\big({\cal R}_2h\big)-{\cal R}_3I\big)\nu=Ae_y\big({\cal R}_1h\big)\nu,\ {\cal P}{\cal R}_2h=0\quad \hbox{on }\partial O,\ea\right.\eeq
whose existence and uniqueness also follows from Lax-Milgram's theorem.
\par
Using these operators, we get that $(\zeta,\eta,q)$ is a solution of (\ref{pbvame}) if and only
taking $h(t)={\cal Q}\zeta(t,.)_{|O}\in {\cal H},$ we have
$$\zeta(t,.)={\cal R}_1h(t),\quad \eta(t,.)={\cal R}_2h(t),\quad q(t)={\cal R}_3h(t),\quad
{dh\over dt}=({\cal R}_2h)_{|O},\quad h(0)=\zeta_0.$$
Therefore, defining ${\cal T}:{\cal H}\to {\cal H}$ by
\beq\label{definTa} {\cal T}h=({\cal R}_2h)_{|O},\eeq
the problem is the existence of solution of the differential equation in ${\cal H}$
$${dh\over dt}={\cal T}h,\quad h(0)=\zeta_0.$$
This just follows from  ${\cal T}$ being a linear and continuous operator in ${\cal H}$. Thus, it defines a  uniformly continuous semigroup ${\cal S}(t)$ of bounded linear operators in ${\cal H}$ by
\beq\label{defseml}{\cal S}(t)h=e^{t{\cal T}}h=\sum_{n=0}^\infty {t^n\over n!}{\cal T}^nh\ \hbox{ in }{\cal H},\quad\forall\, t\geq 0.\eeq
In order to prove (\ref{ideeL1}) we just
use $\partial_t\zeta$ as test function in the first equation in  (\ref{pbvame}) and $\eta$ as test function in the second one. \par
Let us now prove (\ref{decexpl}). For this purpose, we define $(\zeta^l,q^l)\in C^\infty([0,\infty);H^1(O))^N\cap C^\infty([0,\infty);L^2(O)/\RR)^N$ as the solution of the Stokes problem 
$$\left\{\ba{l}\dis -2\mu\, {\rm div}_y e_y(\zeta^l)+\nabla_y q^l=0\ \hbox{ in }(0,\infty)\times O\\ \ecart\dis
{\rm div}_y\zeta^l=0\ \hbox{ in }(0,\infty)\times O\\ \ecart\dis
\zeta^l={\cal Q}\zeta\ \hbox{ on }(0,\infty)\times \partial O,\ea\right.$$
which satisfies
$$\int_O|e_y(\zeta^l)|^2dy\leq C\int_{Y\setminus\overline O}|e_y(\zeta)|^2dy\ \hbox{  in }(0,\infty).$$
for some $C$ independent on $\zeta$. Using this inequality, taking $\zeta$ as test function in the first equation in (\ref{pbvame}),  $\zeta^l$ as test function in the second one and adding the corresponding equalities we easily get
$$\int_{Y\setminus \overline O}Ae_y(\zeta):e_y(\zeta)\,dy\leq C\int_O|e_y(\eta)|^2dy\ \hbox{ in }(0,\infty).$$
Replaced in (\ref{ideeL1}), this proves 
$${d\over dt}\int_{Y\setminus \overline O}Ae_y(\zeta):e_y(\zeta)\,dy\leq -
{2\mu\over C}\int_{Y\setminus \overline O}Ae_y(\zeta):e_y(\zeta)dy\ \hbox{ in }(0,\infty),$$
and then by Gronwall's inequality
\beq\label{desgro1}\int_{Y\setminus \overline O}Ae_y(\zeta):e_y(\zeta)dy\leq e^{-{2\mu\over C}t}
\int_{Y\setminus \overline O}Ae_y(\zeta(0,y)):e_y(\zeta(0,y))dy \hbox{ in }(0,\infty).\eeq
Since $\zeta(0,y)={\cal R}_1\zeta_0$, we get
\beq\label{desgro2}\int_{Y\setminus \overline O}Ae_y(\zeta(0,y)):e_y(\zeta(0,y))dy\leq C\|\zeta_0\|^2_{H^{1\over 2}(\partial 0)^N}.\eeq
On the other hand, using that $(\eta,q)$ satifies
$$\left\{\ba{l}\dis 
-2\mu\,{\rm div}_ye_y(\eta)+\nabla q=0\ \hbox{ in }O
\\ \ecart\dis {\rm div}_y\eta=0\ \hbox{ in }O\\ \ecart\dis
 {\cal P}\eta=0,\  \big(2\mu e_y(\eta)-q I\big)\nu=Ae_y(\zeta)\nu\ \hbox{ on }\partial O,
\ea\right.\qquad \hbox{ in }(0,\infty),$$
we also have 
$$\|\eta\|_{H^1(O)^N}+\|q\|_{L^2(O)}\leq C\|\zeta\|_{H^1(Y\setminus\overline O)^N}\ \hbox{ in }(0,\infty),$$
which combined with  (\ref{desgro1}) and  (\ref{desgro2}) proves (\ref{decexpl}).
\cqfd
\begin{Lem}\label{leexspme2} For every
$E\in L^\infty(0,T;\RR^{N\times N}_s)$,
there exits a unique solution $(\zeta,\eta,q)$ of 
\beq\label{pbvame2}\left\{\ba{l}\dis \zeta\in L^\infty(0,T; H^1_\sharp(Y\setminus\overline O))^N,\ \ 
\eta\in L^\infty(0,T; H^1(O))^N,\ \ q\in L^\infty(0,T;L^2(O))\\ \ecart\dis 
-{\rm div}_y\big(Ae_y(\zeta)\big)=0\ \hbox{ in }(0,T)\times (Y\setminus\overline O)\\ \ecart\dis
-2\mu\,{\rm div}_ye_y(\eta)+\nabla_y q=0\ \hbox{ in }(0,T)\times O
\\ \ecart\dis {\rm div}_y\eta=0\ \hbox{ in }(0,T)\times O\\ \ecart\dis
\zeta\hbox{ periodic in }Y,\  A\big(E+e_y(\zeta)\big)\nu\hbox{ anti-periodic in }Y
\\ \ecart\dis 
\partial_t{\cal Q}\big(Ey+\zeta\big)=\eta\ \hbox{ on }(0,T)\times \partial O\\ \ecart\dis
{\cal P}\big(A\big(E+e_y(\zeta)\big)\nu\big)=0\ \hbox{ on }(0,T)\times \partial O\\ \ecart\dis
A\big(E+e_y(\zeta)\big)\nu=(2\mu\, e_y(\eta)-qI)\nu\ \hbox{ on }(0,T)\times\partial O\\ \ecart\dis
{\cal Q}(Ey+\zeta)_{|t=0}=0\ \hbox{ on }\partial O.
\ea\right.
\eeq
Moreover, there exist a positive and symmetric tensor $A_h\in {\cal L}(\RR^N_s;\RR^N_s)$ and a tensor function $S\in C^\infty([0,\infty);{\cal L}(\RR^N_s;\RR^N_s))$ with
\beq\label{decexS}\|S(t)\|_{{\cal L}(\RR^N_s;\RR^N_s)}\leq \Lambda e^{-\lambda t},\quad\forall\,t\geq 0,\eeq
for some $\lambda,\Lambda>0$, such that
\beq\label{expAhS}\int_{Y\setminus \overline O} A\big(E+e_y(\zeta)\big)dy+\int_O\big(2\mu e_y(\eta)-q I\big)dy=A_hE+\int_0^t S(t-s)E(s)\,ds\ \hbox{ in }(0,T),
\eeq
and, assuming $E$ in $W^{1,1}(0,T;\RR^{N\times N}_s)$
\beq\label{energS}\ba{l}\dis
\left(A_hE+\int_0^t S(t-s)E(s)\,ds\right):{dE\over dt}\\ \ecart\dis ={1\over 2}{d\over dt}\int_{Y\setminus\overline O}A\big(E+e_y(\zeta)\big):\big(E+e_y(\zeta)\big)dy
+2\mu\int_O|e_y(\eta)|^2dy.\ea\eeq
\end{Lem}
\noindent
{\bf Proof.}  Let us first assume that $E$ is in $W^{1,1}(0,T)$. We define $Z:\RR^{N\times N}_s\to H^1_\sharp(Y\setminus \overline O)^N$ by 
\beq\label{deZx}\left\{\ba{l}\dis -{\rm div}_y\big(Ae_y(Z \xi)\big)=0\ \hbox{ in }Y\setminus\overline O\\ \ecart
\dis Z\xi\hbox{ periodic in }Y,\ A\big(\xi+e_y(Z\xi)\big)\nu\hbox{ anti-periodic in }Y\\ \ecart\dis
A\big(\xi+e_y(Z \xi)\big)\nu=0\ \hbox{ on }\partial O,\quad
\int_{Y\setminus\overline O}\hskip-5pt Z\xi\,dy=0,
\ea\right.\qquad\forall\,\xi\in\RR^{N\times N}_s.
\eeq
Then, taking $\tilde \zeta=\zeta-ZE$, we get that (\ref{pbvame2}) is equivalent to
\beq\label{pbvame3}\left\{\ba{l}\dis \tilde\zeta\in W^{1,1}(0,T; H^1_\sharp(Y\setminus\overline O))^N,\ \ 
\eta\in W^{1,1}(0,T; H^1(O))^N,\ \ q\in W^{1,1}(0,T;L^2(O))\\ \ecart\dis 
-{\rm div}_y\big(Ae_y(\tilde\zeta)\big)=0\ \hbox{ in }(0,T)\times (Y\setminus\overline O)\\ \ecart\dis
-2\mu\,{\rm div}_ye_y(\eta)+\nabla_y q=0\ \hbox{ in }(0,T)\times O
\\ \ecart\dis {\rm div}_y\eta=0\ \hbox{ in }(0,T)\times O\\ \ecart\dis
\tilde \zeta\hbox{ periodic in }Y,\  Ae_y(\tilde\zeta)\nu\hbox{ anti-periodic in  }Y
\\ \ecart\dis 
{\cal Q}\partial_t\tilde\zeta=\eta-{\cal Q}\big({dE\over dt}y+Z{dE\over dt}\big)\ \hbox{ on }(0,T)\times \partial O\\ \ecart\dis
{\cal P}\big(Ae_y(\tilde\zeta)\nu\big)=0,\
Ae_y(\tilde \zeta)\nu=(2\mu e_y(\eta)-qI)\nu\ \hbox{ on }(0,T)\times\partial O\\ \ecart\dis
{\cal Q}(Ey+\tilde \zeta+ZE)_{|t=0}=0\ \hbox{ on }\partial O.
\ea\right.
\eeq
Using the operators ${\cal R}$ and ${\cal T}$ defined in the proof of Lemma \ref{leexspme}, this means that $h:={\cal Q}\tilde \zeta_{|\partial O}$ satisfies
\beq\label{ecu2dh} {dh\over dt}={\cal T}h-{\cal Q}\big({dE\over dt}y+Z{dE\over dt}\big)\,\quad h(0)=-{\cal Q}(E(0)y+(ZE)(0,y)).\eeq
and then
$$\tilde \zeta(t,.)={\cal R}_1h(t),\quad \eta(t,.)={\cal R}_2h(t),\quad  q(t,.)={\cal R}_3h(t).$$
Thus, the problem is the existence of solution for (\ref{ecu2dh}). Since ${\cal T}$ defines the semigroup ${\cal S}(t)$, this follows from the semigroup theory, which gives
$$h(t)=-{\cal S}(t){\cal Q}(E(0)y+ZE(0,y))-
\int_0^t{\cal S}(t-s)\partial_t{\cal Q}\big(Ey+ZE\big)\,ds\ \hbox{ in }H^{1\over 2}(\partial O)^N,\quad t\in [0,T],$$ 
or integrating by parts
$$h(t)=-{\cal Q}\big(Ey+ZE\big)
-\int_0^t{\cal S}(t-s){\cal T}{\cal Q}\big(Ey+ZE\big)ds\ \hbox{ in }H^{1\over 2}(\partial O)^N,\quad t\in [0,T].$$
Therefore
$$\zeta=ZE-{\cal R}_1{\cal Q}\big(Ey+ZE\big)-\int_0^t{\cal R}_1{\cal S}(t-s){\cal T}{\cal Q}\big(Ey+ZE\big)ds$$
$$\eta=-{\cal R}_2{\cal Q}\big(Ey+ZE\big)-\int_0^t{\cal R}_2{\cal S}(t-s){\cal T}{\cal Q}\big(Ey+ZE\big)ds$$
$$q=-{\cal R}_3{\cal Q}\big(Ey+ZE\big)-\int_0^t{\cal R}_3{\cal S}(t-s){\cal T}{\cal Q}\big(Ey+ZE\big)ds.$$
We observe that these functions are well defined for $E$ just in $L^\infty(0,T)$ and provide the solution of (\ref{pbvame2}). Moreover these expressions imply (\ref{expAhS}) with
$$\ba{ll}\dis A_h\xi&\dis=\int_{Y\setminus \overline O} A\big(\xi+e_y\big(Z\xi-{\cal R}_1{\cal Q}(\xi y+Z\xi)\big)\big)dy\\ \ecart &\dis -\int_O\big(2\mu e_y\big({\cal R}_2{\cal Q}\big(\xi y+Z\xi\big)-{\cal R}_3{\cal Q}\big(\xi y+Z\xi\big)\big)dy,\quad\forall\,\xi\in\RR^{N\times N}_s,\ea$$
$$\ba{ll}\dis S(s)\xi&\dis=-\int_{Y\setminus \overline O} Ae_y\big({\cal R}_1{\cal S}(s){\cal T}{\cal Q}\big(\xi y+Z\xi\big)\big)dy\\ \ecart &\dis -\int_O\big(2\mu e_y\big({\cal R}_2{\cal S}(s){\cal T}{\cal Q}\big(\xi y+Z\xi\big)\big)-{\cal R}_3{\cal S}(s){\cal T}{\cal Q}\big(\xi y+Z\xi\big)\big)dy,\quad\forall\,\xi\in\RR^{N\times N}_s,\ \forall\,s\geq 0.\ea$$
\par
In order to prove that $A_h$ is symmetric and positive, we take two matrices $\xi_1,\xi_2\in \RR^{N\times N}_s$ and we introduce
$\phi_i={\cal R}_1{\cal Q}(\xi_iy+Z\xi_i)$, $\psi_i={\cal R}_2{\cal Q}(\xi_iy+Z\xi_i)$ and 
$\pi_i={\cal R}_3{\cal Q}(\xi_iy+Z\xi_i),$ $i=1,2$ in such way that
$$A_h\xi_i\cdot \xi_j=\left(\int_{Y\setminus\overline O}A\big(\xi_i+e_y\big(Z\xi_i-\phi_i\big)\big)dy-\int_O\big(2\mu e_y(\psi_i)-\pi_i I\big)dy\right):\xi_j,\quad 1\leq i,j\leq 2.$$
We recall that by (\ref{defR1}), (\ref{defR2}) and a.e. in $(0,T)$, the functions  $\phi_i,\psi_i,\pi_i$ are the solutions of
\beq\label{ecRep1}\left\{\ba{l}\dis-{\rm div}_y(Ae_y(\phi_i))=0\ \hbox{ in }Y\setminus\overline O\\ \ecart\dis
-2\mu\,{\rm div}_ye_y(\psi_i)+\nabla \pi_i=0\ \hbox{ in }O\\ \ecart\dis
{\rm div}_y\psi_i=0\ \hbox{ in }O\\ \ecart\dis
\phi_i\hbox{ periodic in }Y,\ Ae_y(\phi_i)\nu=0\hbox{ anti-periodic in }Y\\ \ecart\dis
{\cal Q}(\phi_i-\xi_iy-Z\xi_i)=0,\ {\cal P}(Ae(\phi_i)\nu)=0 \ \hbox{ on }\partial O\\ 
\ecart\dis Ae_y(\phi_i)\nu= \big(2\mu\,e_y(\psi_i)-\pi_i I\big)\nu,\ {\cal P}\psi_i=0\ \hbox{ on }\partial O.\ea\right.\eeq
We have
$$\ba{l}\dis\int_{Y\setminus \overline O}A\big(\xi_i+e_y(Z\xi_i-\phi_i)\big)dy:\xi_j-
\int_O\big(2\mu e_y(\psi_i)-\pi_i I\big)dy:\xi_j\\ \ecart\dis
=\int_{Y\setminus \overline O}A\big(\xi_i+e_y(Z\xi_i-\phi_i)\big):\big(\xi_j+e_y(Z\xi_j-\phi_j)\big)\,dy\\ \ecart\dis-\int_{Y\setminus \overline O}A\big(\xi_i+e_y(Z\xi_i)\big):e_y(Z\xi_j-\phi_j)dy+\int_{Y\setminus \overline O}Ae_y(\phi_i):e_y(Z\xi_j-\phi_j)dy\\ \ecart\dis -
\int_O\big(2\mu e_y(\psi_i)-\pi_i I\big)dy:\xi_j dy,
\ea$$
where by definition (\ref{deZx}) of $Z$
$$\int_{Y\setminus \overline O}A\big(\xi_i+e_y(Z\xi_i)\big):e_y(Z\xi_j-\phi_j)dy=0$$
and by  (\ref{ecRep1})
$$\ba{l}\dis\int_{Y\setminus \overline O}Ae_y(\phi_i):e_y(Z\xi_j-\phi_j)dy
-\int_O\big(2\mu e_y(\psi_i)-\pi_i I\big)dy:\xi_j\\ \ecart\dis=-\int_{\partial O} Ae_y(\phi_i)\nu\cdot (Z\xi_j-\phi_j)\,d\sigma(y)-\int_{\partial O}\big(2\mu e_y(\psi_i)-\pi_i I\big)\nu\cdot \xi_jy\,d\sigma(y)=0.\ea$$
Therefore, we have proved
\beq\label{posiTD}
\ba{l}\dis\int_{Y\setminus \overline O}A\big(\xi_i+e_y(Z\xi_i-\phi_i)\big):\xi_jdy-
\int_O\big(2\mu e_y(\psi_i)-\pi_i I\big)dy:\xi_j\\ \ecart\dis
=\int_{Y\setminus \overline O}A\big(\xi_i+e_y(Z\xi_i-\phi_i)\big):\big(\xi_j+e_y(Z\xi_j-\phi_j)\big)\,dy,\quad 1\leq i,j\leq 2.
\ea\eeq
This proves that $A_h$ is symmetric and non-negative. Moreover, since 
$Z\xi_i-\phi_i$ periodic implies that it cannot be a linear function except if $\xi_i$ is the null function,  we get that $A_h$ is positive.\par\medskip
Let us now prove (\ref{energS}). We use
\beq\label{posiTDa}\ba{l}\dis\int_{Y\setminus\overline O} A(E+e_y(\zeta))dy :{dE\over dt}=\int_{Y\setminus\overline O} A(E+e_y(\zeta))dy:\Big({dE\over dt}+e_y(\partial_t\zeta)\Big)dy\\ \ecart\dis
-\int_{Y\setminus\overline O} A(E+e_y(\zeta)):e_y(\partial_t\zeta)dy.,\ea\eeq
where using (\ref{pbvame2}), we have
$$\ba{l}\dis \int_{Y\setminus\overline O} A(E+e_y(\zeta)):e_y(\partial_t\zeta)dy=-\int_{\partial O}A(E+e_y(\zeta))\nu:\partial_t\zeta\,d\sigma(y)\\ \ecart\dis =\int_{\partial O}(2\mu \,e_y(\eta)-qI)\nu:\big({dE\over dt} y-\eta\big)\,d\sigma(y=\int_O \big(2\mu\,e_y(\eta)-qI\big):\big({dE\over dt}-e_y(\eta)\big)dy\\ \ecart\dis=\int_O \big(2\mu\,e_y(\eta)-qI\big)dy:{dE\over dt}-2\mu\int_O |e_y(\eta)|^2dy.
\ea$$
Replacing this equality in (\ref{posiTDa}) and using (\ref{expAhS}) we get  (\ref{energS}). \cqfd
\par\medskip
As a consequence of Lemmas \ref{leexspme} and \ref{leexspme2} we get the following macroscopic equation for the function $u$ in Theorem \ref{Thhocde}. It proves that the material obtained by introducing liquid inclusions in an elastic medium behaves as a viscoelastic material.
\begin{Thm} \label{Threpvel} In the conditions of Theorem \ref{Thhocde} and taking $S$ as the function  defined in Lemma \ref{leexspme2},  there exists $R\in C^\infty([0,\infty);{\cal L}(H^{1\over 2}(\partial O)^N;\RR^{N\times N}_s))$ such that the function $u$ in Theorem \ref{Thhocde} satisfies
\beq\label{pblimMs}\left\{\ba{l}\dis \rho^{eff}\partial^2_{tt}u
-{\rm div}_x\Bigg(A_he_x(u)+\int_0^t S(t-s)e_x(u)(s,x)ds\\ \ecart\dis \qquad\qquad+\int_{\partial O}R(t,y)(e_x(u^0)y+u_1^0)d\sigma(y)\Bigg)=f\hbox{ in }(0,T)\times\Om\\ \ecart\dis
u=0\ \hbox{ on }(0,T)\times\partial\Om\\ \ecart\dis
u_{|t=0}=u^0,\quad \rho^{eff}\partial_t u_{|t=0}=\rho^s(1-|O|)u^1+\rho^l|O|v^0.\ea\right.\eeq
Moreover, there exist $\lambda,\Lambda>0$ such that
\beq\label{decexp} \|R(t,.)\|_{{\cal L}(H^{1\over 2}(\partial O)^N;\RR^{N\times N}_s)}\leq \Lambda e^{-\lambda t},\quad \forall\, t\geq 0.\eeq
\end{Thm}
\medskip\noindent
{\bf Proof.} We decompose the functions $(u_1,v,p)$ in (\ref{pbli2e}) as $(u_1,v,p)=(\mathbf{U}_1,\mathbf{V},\mathbf{P})+(\mathsf{U}_1,\mathsf{V},\mathsf{P)}$, where a.e. in $\Om$, we have
$$\left\{\ba{l}\dis 
-{\rm div}_y\big(Ae_y(\mathbf{U}_1)\big)=0\ \hbox{ in }(0,T)\times \big(Y\setminus \overline O\big)\\ \ecart\dis
-2\mu\,{\rm div}_y e_y(\mathbf{V})+\nabla_y\mathbf{P}=0\  \hbox{  in } (0,T)\times O\\ \ecart\dis
{\rm div}_y \mathbf{V}=0\ \hbox{ in } (0,T)\times O\\ \ecart\dis
\mathbf{U}_1\hbox{ periodic in }Y,\ Ae_y(\mathbf{U}_1\big)\nu\hbox{ anti-periodic in }Y\\ \ecart\dis
{\cal Q}\partial_t\mathbf{U}_1=\mathbf{V}\ \hbox{ on }(0,T)\times\partial O\\ \ecart\dis
{\cal P}\big(Ae_y(\mathbf{U}_1\big)\nu\big)=0,\ 
Ae_y(\mathbf{U}_1)\nu=(2\mu e_y(\mathbf{V})-\mathbf{P}I)\nu\ \ \hbox{ on }(0,T)\times \partial O\\ \ecart\dis
({\cal Q}\mathbf{U}_1)_{|t=0}={\cal Q}(e_x(u^0)y+u_1^0)\ \hbox{ on }\partial O,
\ea\right.$$
$$\left\{\ba{l}\dis 
-{\rm div}_y\big(Ae_y(\mathsf{U}_1)\big)=0\ \hbox{ in }(0,T)\times\big(Y\setminus \overline O\big)\\ \ecart\dis
-2\mu\,{\rm div}_y e_y(\mathsf{V})+\nabla_y\mathsf{P}=0\ \hbox{ in }(0,T)\times O\\ \ecart\dis
{\rm div}_y \mathsf{V}=0\ \hbox{ in }(0,T)\times O\\ \ecart\dis
\mathsf{U}_1\hbox{ periodic in }Y,\ A\big(e_x(u)+e_y(\mathsf{U}_1)\big)\nu\hbox{ anti-periodic in }Y\\ \ecart\dis
{\cal Q}\big(e_x(\partial_tu)y+\partial_t\mathsf{U}_1\big)=\mathsf{V}\ \hbox{ on }(0,T)\times\partial O\\ \ecart\dis
{\cal P}\big(A\big(e_x(u)+e_y(\mathsf{U}_1)\big)\nu\big)=0, \ A\big(e_x(u)+e_y(\mathsf{U}_1)\big)\nu=(2\mu e_y(\mathsf{V})-\mathsf{P}I)\nu\ \ \hbox{ on }(0,T)\times\partial O\\ \ecart\dis
{\cal Q}(e_x(u^0)y+\mathsf{U}_1)_{|t=0}=0\ \hbox{ on }\partial O,
\ea\right.$$
By Lemma  \ref{leexspme},  the functions $\mathbf{U}_1$, $\mathbf{V}$, $\mathbf{P}$ are given by
$$\mathbf{U}_1={\cal R}_1{\cal S}(t){\cal Q}\big(e_x(u^0)y+u_1^0\big),\quad \mathbf{V}={\cal R}_2{\cal S}(t){\cal Q}\big(e_x(u^0)y+u_1^0\big),\quad \mathbf{P}={\cal R}_3{\cal S}(t){\cal Q}\big(e_x(u^0)y+u_1^0\big),$$
and therefore, defining $R\in C^\infty([0,\infty);{\cal L}(H^{1\over 2}(\partial O)^N;\RR^{N\times N}_s))$ by 
\beq\label{defiR}\int_{\partial O}R(t,y)\varphi(y)\,d\sigma(y)=\int_{Y\setminus \overline O} Ae_y({\cal R}_1{\cal S}(t){\cal Q}\varphi)\,dy+\int_O\big(2\mu e_y({\cal R}_2{\cal S}(t){\cal Q}\varphi)-{\cal R}_3{\cal S}(t){\cal Q}\varphi\big)dy,\eeq
for every $\varphi\in H^{1\over 2}(\partial O)^N$ and a.e. $t\in (0,T)$, we have
$$
\int_{Y\setminus \overline O} Ae_y(\mathbf{U}_1)dy+\int_O\big(2\mu e_y(\mathbf{V})-\mathbf{P} I\big)dy
=\int_{Y\setminus \overline O}R(t)(e_x(u^0)y+u_1^0)d\sigma(y),
$$
while by Lemma \ref{leexspme2} 
$$
\int_{Y\setminus \overline O} A\big(e_x(u)+e_y(\mathsf{U}_1)\big)dy+\int_O\big(2\mu e_y(\mathsf{V})-\mathsf{P} I\big)dy
=A_he_x(u)+\int_0^tS(t-s)e_x(u)(s)\,ds.
$$
Then, the first equation in (\ref{pbli2e}) gives (\ref{pblimMs}).\par 
Inequality (\ref{decexp}) is a simple consequence of (\ref{decexpl}). \cqfd\par\medskip
\begin{Rem}\label{reinc} The differential equation in (\ref{pblimMs}) contains the term
$$-{\rm div}_x\left(\int_{\partial O}R(t,y)(e_x(u^0)y+u_1^0)d\sigma(y)\right)$$
which is a memory term depending on the initial conditions and acts in the equation as a force depending of such conditions. By the exponential decay of $R$ it has influence only for small times. If we assume that the initial conditions $u_\ep^0$ converge strongly to $u^0$ in $H^1_0(\Om)^N$, then the function $u_1^0$ vanishes and the  equation in (\ref{pblimMs}) reads as
$$\rho^{eff}\partial^2_{tt}u
-{\rm div}_x\Bigg(A_he(u)+\int_0^t S(t-s)e(u)(s,x)ds+\mathsf{R}(t)e_x(u^0)\Bigg)=f\hbox{ in }(0,T)\times\Om,$$
with $\mathsf{R}\in C^\infty([0,\infty);{\cal L}(\RR^{N\times N}_s))$ defined by
$$\mathsf{R}(t)\xi=\int_{\partial O}R(t,y)(\xi y)d\sigma(y),\quad\forall\,\xi\in\RR^{N\times N}_s,\ \forall\,t\in[0,\infty).$$
Another interesting situation corresponds to ${\cal Q}(e_x(u^0)y+u_1^0)=0$ on $\Om\times\partial O$. This is the case if  $u^1_0=-{\cal R}_1{\cal Q}(e_x(u^0)y)$ which, related to the results in \cite{FrMu} (see also \cite{BrCa2}, \cite{CCMM}), we call well prepared initial conditions. In this case the equation in (\ref{pblimMs}) reduces to
$$\rho^{eff}\partial^2_{tt}u
-{\rm div}_x\Bigg(A_he(u)+\int_0^t S(t-s)e(u)(s,x)ds\Bigg)=f\hbox{ in }(0,T)\times\Om.$$
\end{Rem}
As we said in the introduction, our aim in the paper is not just to prove that the limit problem of (\ref{sistho}) contains a memory term such as it is proved in Theorem \ref{Threpvel}. We are also interested in studying the main properties of this nonlocal operator, which are the ones that we could impose on the models that describe the viscoelastic  behavior of skin. The two main  properties in this sense are given by the exponential decay (\ref{decexS}) of $S$ and the positivity condition (\ref{energS}). Let us characterize this last  property in Theorem \ref{posONL} using the Laplace transform. Indeed, we remark that (\ref{energS}) is not
very interesting because it is written using the functions $\zeta$, $\eta$  and $q$ depending on the microstructure. In Theorem \ref{posONL}  we will get the corresponding property just using the macroscopic variable $E$ (see (\ref{prASnl})). We will need to assume that there are no $\xi\in\CC^{N\times N}_s$,  $\xi\not =0$, $k\in\CC$ such that
the overdetermined problem
\beq\label{Dela10a}\left\{\ba{l}\dis 
-{\rm div}_y\big(A\big( \xi+e_y(w_\xi)\big)\big)=0\ \hbox{ in } Y\setminus\overline O\\ \ecart\dis
w_\xi\hbox{ periodic in }Y,\ A\big(\xi+e_y(w_\xi)\big)\nu\hbox{ anti-periodic in }Y
\\ \ecart\dis 
{\cal Q}\big(\xi y+w_\xi\big)=0,\ {\cal P}\big(A\big(\xi+e_y(w_\xi)\big)\nu\big)=0\ \hbox{ on } \partial O\\ \ecart\dis
A\big(\xi+e_y(w_\xi)\big)\nu=-k\nu\ \hbox{ on } \partial O,
\ea\right.
\eeq
has solution. In Proposition \ref{hipsa} below, we will show that this is the case if $O$ is a ball and $A$ is isotropic.
\begin{Pro} \label{posONL} The  Laplace transform $LS(z)$ of the function $S$ defined by Lemma \ref{leexspme2}
is symmetric. \par
If there are no $\xi\not=0$, $k\in\CC$, such that (\ref{Dela10a}) has a solution, then, 
there exists $c>0$ such that
\beq\label{proLS1} Re\Big((A_h+LS(z))\xi:\overline {z\xi}\Big)\geq c\Big(Re(z)+{|z|^2\over 1+|z|^2}\Big)|\xi|^2,\quad\forall\, z\in\CC^N,\ Re(z)>0, \eeq
\beq\label{proLS2} Im\Big((A_h+LS(z))\xi:\overline {z\xi}\Big)\,Im(z) \leq -c\,|Im(z)|^2|\xi|^2,\quad\forall\, z\in\CC^N,\ Re(z)>0. \eeq
Moreover, (\ref{proLS1}) implies that for another constant $c>0$ the tensor functions $A_h$ and $S$ satisfy 
\beq\label{prASnl} \ba{l}\dis\int_0^T \Big(A_hE+\int_0^tS(t-s)E(s)ds\Big):{dE\over dt}\,dt\\ \ecart\dis\geq c|E(T)|^2+c\int_0^T \Big|\int_0^te^{-(t-s)}{dE\over ds}ds\Big|^2dt+c\Big|\int_0^Te^{-(T-s)}{dE\over ds}ds\Big|^2,\ea\eeq
for every $E\in W^{1,1}(0,T;\RR^{N\times N}_s),$ with $E(0)=0$.\end{Pro}\par\noindent
{\bf Proof.} Taking into account (\ref{pbvame2}), (\ref{expAhS}) we deduce that for every $\xi\in \CC^{N\times N}_s$, we have
\beq\label{defLS1}\big(A_h +LS(z)\big)\xi=\int_{Y\setminus \overline O}A\big(\xi+e_y(\hat \zeta)\big)dy+\int_O \big(2\mu e_y(\hat \eta)-\hat q I\big)dy,\quad \forall\, z\in\CC^N,\ Re(z)>0,\eeq
with $(\hat\zeta(z),\hat\eta(z),\hat q(z))$ the solution of
\beq\label{Dela1}\left\{\ba{l}\dis 
-{\rm div}_y\big(Ae_y(\hat\zeta)\big)=0\ \hbox{ in } (Y\setminus\overline O)\\ \ecart\dis
-2\mu\,{\rm div}_ye_y(\hat\eta)+\nabla_y \hat q=0\ \hbox{ in } O
\\ \ecart\dis {\rm div}_y\hat \eta=0\ \hbox{ in }O\\ \ecart\dis
\hat \zeta\hbox{ periodic in }Y,\ A\big(\xi+e_y(\hat\zeta)\big)\nu\ \hbox{ anti-periodic in }Y
\\ \ecart\dis 
z{\cal Q}\big(\xi y+\hat\zeta\big)=\hat\eta\ \hbox{ on } \partial O\\ \ecart\dis
{\cal P}\big(A\big(\xi+e_y(\hat\zeta)\big)\nu\big)=0\ \hbox{ on } \partial O\\ \ecart\dis
A\big(\xi+e_y(\hat\zeta)\big)\nu=(2\mu\, e_y(\hat\eta)-\hat qI)\nu\ \hbox{ on }\partial O.
\ea\right.
\eeq
\par
Taking $\xi_1,\xi_2\in\CC^{N\times N}_s$ and $(\hat\zeta_1,\hat\eta_1,\hat q_1)$, $(\hat\zeta_2,\hat\eta_2,\hat q_2)$ the solutions of (\ref{Dela1}) with $\xi$ replaced by $\xi_1$ and $\xi_2$ respectively, we get
$$\ba{l}\dis z\big(A_h +LS(z)\big)\xi_1:\xi_2 =z\int_{Y\setminus \overline O}A\big(\xi_1+e_y(\hat \zeta_1)\big)dy:\xi_2+z\int_O \big(2\mu e_y(\hat \eta_1)-\hat q_1 I\big)dy:\xi_2\\ \ecart\dis
=z\int_{Y\setminus \overline O}A\big(\xi_1+e_y(\hat \zeta_1)\big):(\xi_2+e_y(\hat\zeta_2)\,dy+z\int_O \big(2\mu e_y(\hat \eta_1)-\hat q_1 I\big)dy:\xi_2\\ \ecart\dis -z\int_{Y\setminus \overline O}A\big(\xi_1+e_y(\hat \zeta_1)\big):e_y(\hat\zeta_2)\,dy,\ea$$
where an integration by parts using (\ref{Dela1}) proves
$$\ba{l}\dis z\int_{Y\setminus \overline O}A\big(\xi_1+e_y(\hat \zeta_1)\big):e_y(\hat\zeta_2)\,dy=-\int_{\partial O} zA\big(\xi_1+e_y(\hat \zeta_1)\big)\nu\cdot\hat \zeta_2\,d\sigma(y)\\ \ecart\dis= -\int_{\partial O} \big(2\mu e_y(\hat\eta_1)-\hat q_1 I\big)\nu\cdot\big(\hat \eta_2-z\xi_2 y\big)d\sigma(y)=-2\mu\int_Oe_y(\hat\eta_1):e_y(\hat\eta_2)\,dy\\ \ecart\dis
+z\int_O\big(2\mu e_y(\hat \eta_1)-\hat q_1I\big)dy:\xi_2.
\ea$$
Thus,
\beq\label{iguALS} z\big(A_h +LS(z)\big)\xi_1:\xi_2 =
z\int_{Y\setminus \overline O}A\big(\xi_1+e_y(\hat \zeta_1)\big)dy:\big(\xi_2+e_y(\hat\zeta_2)\big)dy+2\mu\int_O  e_y(\hat \eta_1):e_y(\hat\eta_2)dy.\eeq
Since $A_h$ is symmetric, this proves that $LS(z)$ is symmetric.\par\medskip
Analogously to (\ref{iguALS}), we also have
 \beq\label{energSLa}
\big(A_h +LS(z)\big)\xi:\overline{z\xi} =\overline z\int_{Y\setminus\overline O}A\big(\xi+e_y(\hat\zeta)\big):\big(\overline{\xi+e_y(\hat\zeta)}\big)dy
+2\mu\int_O|e_y(\hat \eta)|^2dy,\eeq
for every $z\in\CC$, with $Re(z)>0$, and every $\xi\in\CC^{N\times N}_s$.\par
Let us prove the existence of $C>0$ independent of $\xi$ such that
\beq\label{Dela2} |\xi|^2\leq C\int_{Y\setminus \overline O} A\big(\xi+e_y(\hat\zeta)\big):\big(\overline{\xi+e_y(\hat\zeta)}\big)dy\eeq
\beq\label{Dela3}{|z|^2\over 1+|z|^2} |\xi|^2\leq C\int_O|e_y(\hat\eta)|^2dy,\eeq
Defining $Z\xi$ by (\ref{deZx}) (which now takes complex values) we have
\beq\label{Dela4}\ba{l}\dis \int_{Y\setminus\overline O} A(\xi+e_y(Z\xi)):\overline{(\xi+e_y(Z\xi))}dy=\min_{w\in H^1_\sharp(Y\setminus \overline O)}\int_{Y\setminus\overline O} A(\xi+e_y(w))(\overline{\xi+e_y(w)})dy\\ \ecart\dis\leq 
\int_{Y\setminus\overline O}A(\xi+e_y(\hat\zeta)):(\overline{\xi+e_y(\hat \zeta)})dy.\ea\eeq
Since $\xi=e_y(\xi y)$ and $\xi y$ does not belong to $H^1_\sharp(Y\setminus\overline O)$, we get the existence of $C>0$ independent of $\xi$ such that
$$|\xi|^2\leq C\int_{Y\setminus\overline O} A(\xi+e_y(Z\xi)):(\overline{\xi+e_y(Z\xi)})dy.$$
Combined with (\ref{Dela4}), this proves  (\ref{Dela2}).\par
Let us now prove (\ref{Dela3}). It is enough to show that for every $\xi\in\CC^{N\times N}_s$ with $|\xi|=1$, the solution of (\ref{Dela1}) satisfies
\beq\label{Dela5} 0<\inf\left\{\Big(1+{1\over |z|^2}\Big)\int_O|e_y(\hat\eta)|^2dy:\ \ Re(z)>0\right\}.\eeq
We reason by contradiction. If (\ref{Dela5}) does not hold then there exist $z_n$, with $Re(z_n)>0$ and $\xi_n\in \CC^{N\times N}_s$, with $|\xi_n|=1$ such that the solution $\big(\hat\zeta_n, \hat\eta_n,\hat q_n)$ of (\ref{Dela1}) with $\xi$ replaced by $\xi_n$ and $z$ replaced by $z_n$ satisfies
\beq\label{Dela8} \Big(1+{1\over |z_n|}\Big)\hat\eta_n\to 0\ \hbox{ in }H^1(O)^N.\eeq 
Clearly, we can also assume that $\xi_n$ tends to a symmetric matrix function $\xi\in\CC^N$ with $|\xi|=1$. This easily allows us to deduce that for a subsequence (distinguish the case where $z_n$ tends to zero and the one that does not)
\beq\label{Dela9} \hat\zeta_n\to w_\xi\ \hbox{ in }H^1_\sharp(O)^N,\quad \hat q_n\to k\eeq 
with $w_\xi\in H^1_\sharp(Y\setminus \overline O)$, $k\in \CC$ a solution of (\ref{Dela10a}). This contradicts our hypothesis and therefore proves (\ref{Dela3}).\par
Taking the real and imaginary parts in (\ref{energSLa}) and using (\ref{Dela2}) and (\ref{Dela3}) we conclude (\ref{proLS1}) and (\ref{proLS2}).\par
\medskip
Let us now prove (\ref{prASnl}) using (\ref{proLS1}). For this purpose, given $E\in W^{1,1}(0,T;\RR^{N\times N}_s)$ with $E(0)=0$ and $T>0$, we define $E_T:(0,\infty)\to \RR^{N\times N}_s$ by
$$E_T(t)=\left\{\ba{ll}\dis E(t) &\hbox{ if }0\leq t\leq T\\ \ecart\dis E(T) &\hbox{ if }t\geq T.\ea\right.$$ 
Then, recalling that for a function $\varphi$ such that $e^{-at}\varphi\in L^1(0,\infty)$ we have
$$(L\varphi)(a+ib)={\cal F}(e^{-at}\varphi\chi_{(0,\infty)})(b),\quad \forall\, b\in\RR,$$
with ${\cal F}$ the Fourier transform, and that the Fourier transform is an isomorphism in $L^2(\CC)$, we have for every $a>0$ 
\beq\label{Dela10}\ba{l}\dis\int_0^T e^{-2at} \Big(A_hE+\int_0^tS(t-s)E(s)ds\Big):{dE\over dt}\,dt\\ \ecart\dis=\int_\RR 
\big(A_h+LS(a+ib))LE_T(a+ib)\cdot \overline {(a+ib)\,LE_T(a+ib)}\,db
\\ \ecart\dis\geq c\int_\RR \Big(a+{a^2+b^2\over 1+a^2+b^2}\Big)\big|LE_T(a+ib)\big|^2db.
\ea\eeq
To estimate the right-hand side of (\ref{Dela10}) we observe that
$${a^2+b^2\over 1+a^2+b^2}\big|LE_T(a+ib)\big|^2\geq {a^2+b^2\over (1+a)^2+b^2}\big| LE_T(a+ib)\big|^2= \Big|{a+ib\over 1+a+ib}LE_T(a+ib)\Big|^2.
$$
where 
$${a+ib\over 1+a+ib}LE_T(a+ib)=LV(a+ib)$$
with $V$ the solution of
$${dV\over dt}+V={dE_T\over dt},\quad V(0)=0,$$
i.e.
$$V(t)=\int_0^te^{-(t-s)}{dE_T\over ds}\,ds.$$
Returning to (\ref{Dela10}), we then have
\beq\label{Dela11}\ba{l}\dis\int_0^T e^{-2at} \Big(A_hE+\int_0^tS(t-s)E(S)ds\Big){dE\over dt}\,dt\\ \ecart\dis \geq 
c\int_0^\infty e^{-2at}\Big(a|E_T|^2+\Big|\int_0^te^{-(t-s)}{dE_T\over ds}\Big|^2\Big)dt.
\ea\eeq
For the first term in the right-hand side of this inequality, we use
$$\ba{l}\dis a\int_0^\infty e^{-2at}|E_T|^2dt=a\int_0^T e^{-2at}|E(t)|^2dt+a|E(T)|^2\int_T^\infty e^{-2at}dt\\ \ecart\dis = a\int_0^T e^{-2at}|E(t)|^2dt+{e^{-2aT}\over 2}|E(T)|^2.
\ea$$
For the second term, we have
$$\ba{l}\dis \int_0^\infty e^{-2at}\Big|\int_0^te^{-(t-s)}{dE_T\over ds}\Big|^2dt= \int_0^T e^{-2at}\Big|\int_0^te^{-(t-s)}{dE\over ds}\Big|^2dt+\Big|\int_0^Te^{s}{dE\over ds}\Big|^2\int_T^\infty e^{-(2a+1)t}dt\\ \ecart\dis = \int_0^T e^{-2at}\Big|\int_0^te^{-(t-s)}{dE\over ds}ds\Big|^2dt+{e^{-2aT}\over 2a+1}\Big|\int_0^Te^{-(T-s)}{dE\over ds}ds\Big|^2.\ea
$$
Therefore (\ref{Dela11}) implies
$$\ba{l}\dis\int_0^T e^{-2at} \Big(A_hE+\int_0^tS(t-s)E(S)ds\Big){dE\over dt}\,dt\\ \ecart\dis \geq 
ca\int_0^T e^{-2at}|E(t)|^2dt+{ce^{-2aT}\over 2}|E(T)|^2\\ \ecart\dis+c\int_0^T e^{-2at}\Big|\int_0^te^{-(t-s)}{dE\over ds}ds\Big|^2dt+c{e^{-2aT}\over 2a+1}\Big|\int_0^Te^{-(T-s)}{dE\over ds}ds\Big|^2.
\ea$$
Taking $a$ tending to zero in this inequality we conclude (\ref{prASnl}). \cqfd
\begin{Rem} The fact that $LS$ is symmetric is equivalent to have $S$ symmetric.
\end{Rem}

\begin{Pro} \label{hipsa} Assume that $O$ is a ball and that $A$ is isotropic. Then problem (\ref{Dela10a}) has no solution for any $\xi\in\CC^{N\times N}_s$,  $\xi\not= 0$, $k\in\CC$.
\end{Pro}
\par\medskip\noindent
{\bf Proof.} Assume that $O$ agrees with the ball $B_R$ of center $0$ and radius $R\in (0,1/2)$ and that there exists $\lambda^s\geq 0,\mu^s>0$ such that (\ref{isoma}) holds. Reasoning by contradiction, we also assume  that there exists $\xi\in\CC^{N\times N}_s$, $\xi\not=0$, $k\in \CC$, $w_\xi\in H^1_\sharp(Y\setminus\overline{B_R})^N$ solution of (\ref{Dela10a}).  The condition ${\cal Q}(\xi y+w_\xi)=0$ on $\partial B_R$ means that there exists $a\in \CC^N$ and $B\in \CC^{N\times N}$ skew-symmetric, such that $\xi y+w_\xi=a+By$ on $\partial B_R$. Then, the function 
$$\phi(y)=\xi y+w_\xi-a-By,\quad y\in Y\setminus\overline{B_R},$$ 
satisfies the conditions
\beq\label{conphib}-{\rm div}_y\big(Ae_y(\phi)\big)=0 \ \hbox{ in }Y\setminus \overline{B_R},\quad \phi=0, \ Ae_y(\phi)\nu=k\nu\ \hbox{ on }\partial B_R.\eeq
Using polar coordinates, we can construct another solution of (\ref{conphib}) as 
$$\psi(y)=h(|y|)y, $$
with $h$ the unique solution of the Cauchy problem
$$\left\{\ba{l}\dis-{d\over dr}\Big(r^{N+1}(\lambda+2\mu){dh\over dr}+r^N(\lambda N+2\mu)h\Big)\\ \ecart\dis\qquad\qquad
+r^N(\lambda N+2\mu){dh\over dr}+r^{N-1}N(N\lambda+2\mu)h=0\ \hbox{ in }(0,R),\\ \ecart\dis
h(R)=0,\quad (2\mu+\lambda)Rh'(R)+(\lambda N+2\mu)h(R)=k,\ea\right.$$
Thus, $\phi-\psi$ solves
$$-{\rm div}_y\big(Ae_y(\phi-\psi)\big)=0 \ \hbox{ in }Y\setminus \overline {B_R},\quad \phi-\psi=0, \ Ae_y(\phi-\psi)\nu=r\nu\ \hbox{ on }\partial B_R.$$
By unique continuation we conclude that $\phi-\psi$ vanishes on $Y\setminus \overline B_R$ and then 
$$w_\xi=-\xi y+a+By+h(|y|)y\ \hbox{ in }Y\setminus\overline{B_R},$$
but then $w_\xi$ is not periodic in $Y$, which gives a contradiction. \cqfd
\par\medskip
We finish the paper with a simple example of symmetric tensor  $A_h\in {\cal L}(\RR^{N\times N}_s)$ and $S\in C^\infty([0,\infty);{\cal L}(\RR^{N\times N}_s))$ (not necessarily obtained by the homogenization process), such that  conditions (\ref{decexS}), (\ref{proLS1}) and (\ref{proLS2}) are satisfied.
\begin{Pro} \label{Prex} Assume $A_h,B$ two positive symmetric tensors in ${\cal L}(\RR^{N\times N}_s)$ and $\gamma>0$ such that
\beq\label{condal}\gamma A_h-B\ \hbox{ is }positive.\eeq
Then, $A_h$ and $S(t):=-Be^{-\gamma t},$ satisfy (\ref{proLS1}) and (\ref{proLS2}),
with $LS$ the Laplace transform of $S$.
\end{Pro}
\par\medskip\noindent
{\bf Proof.} Taking into account that
$$LS(a+ib)=-{1\over \gamma+a+ib}B=-{\gamma+a-ib\over (\gamma+a)^2+b^2}B,\quad \forall\, a>0,\ b\in\RR,$$
we have
$$(a-ib)(A_h+LS(a+ib))=a A_h-{a(\gamma+a)-b^2\over (\gamma+a)^2+b^2}B-bA_hi+{b\gamma\over (\gamma+a)^2+b^2}Bi$$
and thus
$$Re\big((a-ib)(A_h+LS(a+ib))\big)=a\Big(A_h-{\gamma+2a\over 
(a+\gamma)^2+b^2}B\Big)+{a^2+b^2\over (a+\gamma)^2+b^2}B$$
$$b\,Im\big((a-ib)(A_h+LS(a+ib))\big)=-b^2\Big(A_h-{\gamma\over (\gamma+a)^2+b^2}B\Big)$$
Since
$$\max_{a\geq 0,\,b\in\RR}{\gamma+2a\over (\gamma+a)^2+b^2}={1\over \gamma},$$
and (\ref{condal}) holds, there exists $c>0$ such that
$$\Big(A_h-{\gamma\over (\gamma+a)^2+b^2}B\Big)\xi\cdot\overline \xi\geq\Big(A_h-{\gamma+2a\over (\gamma+a)^2+b^2}B\Big)\xi\cdot\overline \xi\geq c|\xi|^2,\quad \forall\xi\in \CC^{N\times N}_s,$$
which combined with the positivity of $B$ proves (\ref{proLS1}) and (\ref{proLS2}). \cqfd
\begin{Rem} For $A_h, B, \gamma$ In the conditions of Proposition \ref{Prex} and $\rho>0$, we remark that problem
\beq\label{pbej}\left\{\ba{l}\dis\rho\partial^2_{tt}u-{\rm div}_x\Big(A_he_x(u)-\int_0^t Be^{-\gamma(t-s)}\, e_x(u)ds\Big)=f\ \hbox{ in }(0,T)\times \Om\\ \ecart\dis u=0\ \hbox{ on }\partial\Om\\ \ecart\dis u_{|t=0}=u^0,\ (\partial_tu)_{|t=0}=u^1\ \hbox{ in }\Om,\ea\right.\eeq
satisfies an energy conservation law. For this purpose, it is better to integrate by parts the integral term in the equation to get
$$\rho\partial^2_{tt}u-{\rm div}_x\Big(\Big(A_h-{1\over \gamma}B\Big)e_x(u)+{1\over \gamma}\int_0^t Be^{-\gamma(t-s)}\, e_x(\partial_su)ds\Big)=f+{e^{-\gamma t}\over\gamma}{\rm div}_x\big(Be_x(u^0)\big)\ \hbox{ in }(0,T)\times \Om,$$
which taking
$$w=\int_0^te^{-\gamma(t-s)}u\,ds,$$
allows us to write the differential equation as the system
\beq\label{siuw}\left\{\ba{l}\dis \rho\partial^2_{tt}u-{\rm div}_x\Big(\Big(A_h-{1\over \gamma}B\Big)e_x(u)+{1\over \gamma} B e_x(\partial_tw)\Big)=f\ \hbox{ in }(0,T)\times \Om\\ \ecart\dis
\partial_tw+\gamma w=u\ \hbox{ in }(0,T)\times \Om.
\ea\right.\eeq
Multiplying (formally) the first equation by $\partial_t u$ and integrating in $\Om$ we get
$${1\over 2}{d\over dt}\into \Big(\rho |\partial_tu|^2+\Big(A_h-{1\over \gamma}B\Big)e_x(u):e_x(u)\Big)dx+{1\over \gamma}\into Be_x(\partial_t w):e_x(\partial_t u)\,dx=\into f\cdot\partial_t u\,dx,$$
where thanks to the second equation in (\ref{siuw}), we have
$$\into Be_x(\partial_t w):e_x(\partial_t u)\,dx=\into Be_x(\partial_t w):e_x(\partial^2_{tt} w)\,dx+\gamma\into B e_x(\partial_t w):e_x(\partial_tw)\,,dx. $$
Therefore, we get the energy identity
$$\ba{l}\dis{1\over 2}{d\over dt}\into \Big(\rho |\partial_tu|^2+\Big(A_h-{1\over \gamma}B\Big)e_x(u):e_x(u)
+{1\over \gamma}Be_x(\partial_t w):e_x(\partial_t w)
\Big)dx\\ \ecart\dis\qquad\qquad \qquad\qquad\qquad \qquad+\into Be_x(\partial_t w):e_x(\partial_t w)\,dx=\into f\cdot\partial_t u\,dx.\ea$$
Assuming $f\in L^1(0,T;L^2(\Om))^N,$ $u^0\in H^1_0(\Om)^N$, $u^1\in L^2(\Om)^N$, this allows us in particular to deduce  the existence of a unique solution $u$ for (\ref{pbej}) in $W^{1,\infty}(0,T;L^2(\Om))^N\cap L^\infty(0,T; H^1_0(\Om))^N$.
\end{Rem}

\section*{Aknowledgments} This work has been partially supported by the project PID2020-116809GB-I00 of the {\it Ministerio de Ciencia e Innovaci\'on}  of the governement of Spain.


\begin{thebibliography}{20}

\bibitem{AKPS}{\sc Z.~Abdessamad, I.~Kostin, G.~Panasenko, V.P.~Smyshlyaev.} {\it Memory effect in homogenization of a viscoelastic Kelvin-Voigt model with time-dependent coefficients.} Math. Models Methods Appl. Sci. {\bf 19}  (2009), 1603-1630.

\bibitem{All}{\sc G.~Allaire.} {\it Homogenization and two-scale convergence.} SIAM J. Math. Anal. {\bf 23} (1992), 1482-1518.

\bibitem{ArDoHo} {\sc T. Arbogast, J. Douglas, U. Hornung.} {\it Derivation of the double porosity model of single phase flow via homogenization theory.} SIAM J. Math. Anal. {\bf 21} (1990), 823--836.

\bibitem{AGMR}{\sc A.~\'Avila, G.~Griso, B.~Miara, E.~Rohan.} {\it Multiscale modeling of elastic waves: theoretical justification and numerical simulation of band gaps.}  Multiscale Model. Simul. {\bf 7} (2008), 1-21.

\bibitem{BrCa}{\sc M.~Briane, J.~Casado-D\'{\i}az.} {\it Homogenization of an elastodynamic system with a strong magnetic field and soft inclusions inducing a viscoelastic effective behavior.}  J. Math. Anal. Appl. {\bf 492} (2020), 124472, 24 pp.

\bibitem{BrCa2}{\sc M.~Briane, J.~Casado-D\'{\i}az.} {\it Increase of mass and nonlocal effects in the homogenization of magneto-elastodynamics problems.} Calc. Var. Partial Differential Equations {\bf 60} (2021).

\bibitem{Cas} {\sc J. Casado-D\'{\i}az.} {\it Two-scale convergence for nonlinear Dirichlet problems in perforated domains.} Proc. Roy. Soc. Edimburgh A, {\bf 130} (2000), 249--276.

\bibitem{CFLM} {\sc J. Casado-D\'{\i}az, G.A. Francfort, O. L\'opez-Pamies, M.G. Mora.} {\it Liquid filled elastomers: from linearization to elastic enhancement.}  ArXiv:2309.03630v1.

\bibitem{CCMM} {\sc J. Casado-D\'{\i}az, J. Couce-Calvo, F. Maestre, J.D. Mart\'{\i}n-G\'omez.} 
{\it Homogenization and corrector for the wave equation with periodic coefficients.} Math. Models Methods Appl. Sci. {\bf 24} (2014), 1343--1388.

\bibitem{CiDaGr} {\sc D. Cioranescu, A. Damlamian, G. Griso.} {\it Periodic unfolding and homogenization.} C.R. Acad. Sci. Paris, Ser. I {\bf 335} (2002), 99--104.

\bibitem{DMSa} {\sc G. Dal Maso, F. Sapio.} {\it Quasistatic limit of a dinamic viscoelastic model with memory.} 
Milan J. Math. {\bf 89} (2021), 485--522.

\bibitem{FrMu} {\sc G. A. Francfort, F. Murat.}  {\it Oscillations and energy densities in the wave equation.} Commun. Partial Differential Equations {\bf 17} (1992) 1785--1865.
\bibitem{Mur} {\sc J.D. Murray.} {\it Mathematical Biology II: spatial models and biomedical applications.} Interdisciplinary Applied Mathematics {\bf 18}, Springer-Verlag, New York, 2013.

\bibitem{Ngu}{\sc G.~Nguetseng.} {\it A general convergence result for a functional related to the theory of homogenization.} SIAM J. Math. Anal.  {\bf 20} (3) (1989), 608-623.

\bibitem{OlShYo} {\sc O. Oleinik, A.S. Shamaev, G.A. Yosifian.} {\it Mathematical problems in elasticity and homogenization.} Studies in Mathematics and its Applications {\bf 26}. North-Holland Publishing Co. Amsterdam, 1992.

\bibitem{San}{\sc E.~S\'anchez-Palencia.} {\it Nonhomogeneous media and vibration theory}. Lecture Notes in Physics {\bf 127}. Springer-Verlag, Berlin-New York, 1980.


\end{thebibliography}
\end{document}